\documentclass{article}

\usepackage{amssymb,amsmath,mathtools}
\usepackage{authblk}

\usepackage{amsthm}
\usepackage[usenames,dvipsnames,svgnames,table]{xcolor}
\usepackage{algorithm2e}

\newtheorem{theorem}{Theorem}[section]
\newtheorem{proposition}[theorem]{Proposition}

\newtheorem{lemma}[theorem]{Lemma}

\theoremstyle{definition}
\newtheorem{definition}[theorem]{Definition}

\theoremstyle{remark}
\newtheorem{remark}[theorem]{Remark}

\numberwithin{equation}{section}

\newcommand{\argmin}{\operatornamewithlimits{arg\, min}}

\usepackage{enumerate}

\begin{document}

\title{A trust-region SQP method for the numerical approximation of viscoplastic fluid flow}

\author{Timm Treskatis\thanks{Corresponding author. \emph{timm.treskatis@pg.canterbury.ac.nz}}, Miguel A. Moyers-Gonzalez, Chris J. Price}
\affil{School of Mathematics and Statistics, University of Canterbury, Private Bag 4800, Christchurch 8140, New  Zealand}
\date{16 October 2014}
\maketitle
\begin{abstract}
We present a new approach to the problem of stationary viscoplastic duct flow as modelled by the Herschel-Bulkley model, with Bingham fluids included as a special case.

While the mathematical formulation of this problem is conventionally based on a variational inequality, or equivalently, on a nonsmooth minimisation problem for the flow velocity, we suggest an alternative approach. Considering the Lagrangian dual in terms of the stress, rather than the velocity, turns out to be advantageous in numerous ways. The objective functional possesses higher regularity, which ensures applicability of second order methods. Our numerical experiments with a trust-region SQP algorithm also demonstrate clearly superior performance compared to the widely used augmented Lagrangian method, although no artificial regularisation is introduced into the problem.

Hence, besides providing a new theoretical angle to a classical problem, our results also pave the way for an entirely new class of numerical approaches to simulating flows of viscoplastic fluids.
\end{abstract}

\section{Introduction}
\label{sec-intro}

Viscoplastic fluids are characterised by the existence of a yield stress, which marks the transition between viscous and plastic material behaviour. Examples of such materials are encountered in the consumer goods industry (toothpaste, paint), particularly in food processing (dough, tomato sauce) \cite{Bird1983}. Viscoplastic models are also employed in geology and geophysics, for instance to describe flows of lava \cite{Balmforth2000}, lahars \cite{Manville1998} or liquefied soil after a major earthquake \cite{Uzuoka1998}.

The mathematical models named after Bingham \cite{Bingham1922} or Herschel and Bulkley \cite{Herschel1926} lead to problems of free boundary type, as a description of the interfaces between liquid and solid phases is only implicitly contained in the governing equations. The nonsmooth transition between these two regimes and the non-uniqueness of the stress in plastic regions have always posed major obstacles for numerical simulations of viscoplastic fluids.

Approaches to computational solutions of viscoplastic flow problems fall into two distinct categories.

Methods that approximate the viscoplastic models with a purely viscous regularisation were historically first developed. In this context we mention the work of Bercovier and Engelman \cite{Bercovier1980}, the bi-viscosity model of Tanner and Milthorpe \cite{Tanner1983} and the widely used Papanastasiou regularisation \cite{Papanastasiou1987}. More recently, de los Reyes and Gonzalez Andrade \cite{Reyes2009,Reyes2010} employed Fenchel's theory of duality and a Tikhonov regularisation of the dual problem to derive an algorithm that estimates the locations of yielded and unyielded regions as a by-product.

All of these approaches share the advantage of fast convergence, since methods of Newton type are directly applicable. They however become increasingly ill-conditioned the smaller the regularisation is chosen, i.e. the better viscoplasticity is approximated. Consequently, one must seek a compromise between stability and exactness of the method. Furthermore, replacing plasticity with a large, but not infinite viscosity means that methods of this first category are bound to fail with predicting circumstances under which a strictly viscoplastic flow would stop. This inadequacy with reflecting plasticity was demonstrated by Moyers-Gonz{\'a}lez and Frigaard in \cite{Moyers2004}.

Numerical methods from a second class solve the governing equations without prior regularisation. Duvaut and Lions \cite{Duvaut1972,Duvaut1976} established a variational form of the Bingham flow problem in terms of a nonsmooth minimisation problem, or equivalently, a variational inequality. This pioneering work forms the basis of an augmented Lagrangian formulation of the problem and a corresponding numerical strategy commonly named ALG2. Originally, Glowinski \cite{Glowinski1984} proposed this alternating direction method, which is widely used nowadays. For other techniques that avoid regularisation, such as further algorithms of Uzawa type or pseudo time relaxation, we refer to the review of Dean et al. \cite{Dean2007} and the references therein. More recently, Aposporidis et al. \cite{Aposporidis2011} applied Picard iterations to a mixed formulation of the viscoplastic flow problem, which may or may not include regularisation. A few more unregularised algorithms are known \cite{Beris1985,Szabo1992} which, however, are only valid for a very limited range of applications.

Not introducing any artificial regularisation to the problem means that plastic behaviour is exactly represented in the numerical solution, while the computation remains robust and well-conditioned. However, convergence tends to be significantly slower. In particular, methods like ALG2 and classical Uzawa algorithms do not achieve a quadratic or even superlinear rate of convergence. In addition, for methods that apply augmented Lagrangian techniques, the optimal choice of free parameters in the algorithm, which strongly affect the speed of convergence, is still an open problem \cite[pp~126-127]{Fortin1983}.

In our work, we follow a different approach. We first present a dual formulation for stationary Bingham or Herschel-Bulkley flow in terms of a linearly constrained minimisation problem. The objective functional possesses higher regularity compared to the well-known classical functional of Duvaut and Lions. This feature opens up new possibilities for numerical solutions, by bringing second-order methods into play. In this paper, we choose to first discretise the optimisation problem with finite elements and propose a trust-region SQP method to tackle the resulting finite-dimensional problem.

This paper is organised in the following manner: in Section \ref{sec-problem} we introduce the model equations and a corresponding weak formulation. We also provide a summary of key results on the existence and uniqueness of solutions. Once we have presented the new, dual minimisation problem in Section \ref{sec-dual} and its discretisation in Section \ref{sec-discretisation}, Section \ref{sec-algorithm} is devoted to the numerical algorithm for its solution. Finally, we assess the performance of the new approach by carrying out a number of numerical experiments.

\section{Problem statement}
\label{sec-problem}

\subsection{Strong formulation}

We consider the problem of flow through an infinitely long duct with homogeneous cross-section. We assume the cross-sectional domain $\Omega \subset \mathbb{R}^2$ to be Lipschitz with boundary $\partial \Omega = \Gamma$. Furthermore, we let the boundary be decomposed into measurable, disjoint Dirichlet and Neumann sections $\Gamma = \Gamma_\mathrm{D} \mathbin{\dot{\cup}} \Gamma_\mathrm{N}$, where $\Gamma_\mathrm{D}$ is required to possess positive measure.

The Bingham and Herschel-Bulkley models provide constitutive relations between the stress $\boldsymbol\tau$ and the rate of strain $\dot{\boldsymbol\gamma}$; both are vector fields $\Omega \to \mathbb{R}^2$. For the case of duct flow considered here, the rate of strain is given by $\nabla y$, the gradient of the scalar flow velocity $y: \Omega \to \mathbb{R}$ through the cross-section. Model parameters include a plasticity threshold or yield stress $\tau_0 \geq 0$, for Bingham fluids additionally a plastic viscosity parameter $\mu > 0$ and for Herschel-Bulkley fluids a consistency $\kappa > 0$ as well as a power-law exponent $\alpha > 1$.

The flow is driven by a pressure gradient or volumetric force density $f: \Omega \to \mathbb{R}$.

In the classical strong formulation, stationary duct flow of a Herschel-Bulkley fluid is governed by the system
\begin{subequations}
\label{eq-hb-c}
\begin{align}
\boldsymbol\tau &= \kappa \left\vert \dot{\boldsymbol\gamma} \right\vert^{\alpha - 2} \dot{\boldsymbol\gamma} + \tau_0 \frac{\dot{\boldsymbol\gamma}}{\left\vert \dot{\boldsymbol\gamma} \right\vert} & &\mbox{if } \dot{\boldsymbol\gamma} \neq 0 \label{eq-hb-c-visc}\\
\left\vert {\boldsymbol\tau} \right\vert &\leq \tau_0 & &\mbox{if } \dot{\boldsymbol\gamma} = 0 \label{eq-hb-c-plas}\\
\intertext{with conservation of momentum}
-\mathrm{div}\, {\boldsymbol\tau} &= f & &\mbox{in } \Omega \label{eq-hb-c-com}\\
\intertext{and either no slip or perfect slip (symmetry) boundary conditions}
y &= 0 & &\mbox{on } \Gamma_\mathrm{D} \label{eq-hb-c-dbc}\\
\boldsymbol\tau \cdot \vec{\boldsymbol n} &= 0 & &\mbox{on } \Gamma_\mathrm{N}. \label{eq-hb-c-nbc}
\end{align}
\end{subequations}
Here $\vec{\boldsymbol n}$ denotes the outward normal unit vector on the boundary section $\Gamma_\mathrm{N}$.

The special case of a Bingham fluid is recovered through the choices $\alpha = 2$ and $\mu = \kappa$. \eqref{eq-hb-c-visc} and \eqref{eq-hb-c-plas} then assume the simplified form
\begin{subequations}
\label{eq-bi-c}
\begin{align}
\phantom{-\mathrm{div}\,} \boldsymbol\tau &= \mu \dot{\boldsymbol\gamma} + \tau_0 \frac{\dot{\boldsymbol\gamma}}{\left\vert \dot{\boldsymbol\gamma} \right\vert} \phantom{\left\vert \dot{\boldsymbol\gamma} \right\vert^{\alpha - 1}} & &\mbox{if } \dot{\boldsymbol\gamma} \neq 0 \label{eq-bi-c-visc}\\
\left\vert {\boldsymbol\tau} \right\vert &\leq \tau_0 & &\mbox{if } \dot{\boldsymbol\gamma} = 0. \label{eq-bi-c-plas}
\end{align}
\end{subequations}

\subsection{Weak formulation}

A weak formulation of the Bingham flow problem is introduced in \cite{Duvaut1972,Duvaut1976}. An extension to general Herschel-Bulkley fluids is given in \cite{Huilgol2005}. Similar to these works, we consider the following function spaces:

For Lebesgue spaces of order $\alpha$ we apply the standard notation $L^\alpha (\Omega)$. Spaces of $\mathbb{R}^2$-valued functions are set in boldface. Equipped with the canonical norms $\left\Vert \cdot \right\Vert_\alpha$, $L^\alpha (\Omega)$ and $\boldsymbol L^\alpha (\Omega)$ are Banach spaces. If $\alpha = 2$, i.e. within the Bingham scenario, we obtain Hilbert space structure with respect to the $L^2$ or $\boldsymbol L^2$ scalar product $\left( \cdot, \cdot \right)$. The pairing between $L^\alpha (\Omega)$ or $\boldsymbol L^\alpha (\Omega)$ and their dual spaces $L^{\alpha'} (\Omega)$ or $\boldsymbol L^{\alpha'} (\Omega)$, respectively, is referred to as $\left( \cdot, \cdot \right)_{\alpha}$, where the dual index $\alpha'$ is defined through $1/\alpha + 1/\alpha' = 1$. As usual, dual spaces are marked with an asterisk.

With $W^{1,\alpha} (\Omega)$ we denote the Sobolev space of $L^\alpha(\Omega)$-functions with first derivative in $L^\alpha(\Omega)$.

For the Euclidean scalar product in $\mathbb{R}^d$, we use the notation $\cdot$, for the $\alpha$-norm in $\mathbb{R}^d$ we use $\vert \cdot \vert_\alpha$. If $\alpha = 2$, we may omit the subscript. The symbol $C$ is used as a generic constant.

Furthermore, we introduce the subspace of admissible velocity fields and velocity test functions
\begin{align*}
Y &:= \left\lbrace y \in W^{1,\alpha} (\Omega) : y \vert_{\Gamma_\mathrm{D}} = 0 \right\rbrace \subset W^{1,\alpha} (\Omega)
\intertext{For the ease of notation, we define the bi-linear, if $\alpha = 2$, otherwise non-linear form $a : Y \times Y \to \mathbb{R}$}
a(y,z) &:= \kappa \left( \left\vert \nabla y \right\vert^{\alpha - 2} \nabla y, \nabla z \right)_\alpha
\intertext{and the non-linear functionals $j: Y \to \mathbb{R}$}
j(y) &:= \tau_0 \int\limits_\Omega \left\vert \nabla y \right\vert \,\mathrm{d}x
\intertext{as well as $c: Y \to \mathbb{R}$}
c(y) &:= \frac{1}{\alpha} a(y,y) = \frac{\kappa}{\alpha} \int\limits_\Omega \left\vert \nabla y \right\vert^{\alpha} \,\mathrm{d}x.
\end{align*}

We justify that with the above definitions, $a$ and $j$ are indeed well-defined.

\begin{remark}
\begin{enumerate}[(a)]
\item{
Since $y \in W^{1,\alpha} (\Omega)$, we have $\nabla y \in \boldsymbol L^{\alpha} (\Omega)$ and thus
\begin{align*}
\left\Vert \left\vert \nabla y \right\vert^{\alpha - 2} \nabla y \right\Vert_{\alpha'}^{\alpha'} &= \int\limits_\Omega \left\vert \left\vert \nabla y \right\vert_2^{\alpha - 2} \nabla y \right\vert_{\alpha'}^{\alpha'} \,\mathrm{d}x\\
&\leq C \int\limits_\Omega \left\vert \left\vert \nabla y \right\vert_2^{\alpha - 2} \nabla y \right\vert_{2}^{\alpha'} \,\mathrm{d}x\\
&= C \int\limits_\Omega \left\vert \nabla y \right\vert_2^{(\alpha - 2) \alpha'} \left\vert \nabla y \right\vert_2^{\alpha'} \,\mathrm{d}x\\
&= C \int\limits_\Omega \left\vert \nabla y \right\vert_2^{(\alpha - 1) \frac{\alpha}{\alpha - 1}} \,\mathrm{d}x\\
&= C \int\limits_\Omega \left\vert \nabla y \right\vert_2^{\alpha} \,\mathrm{d}x < \infty.
\end{align*}
Consequently, $\left\vert \nabla y \right\vert^{\alpha - 2} \nabla y \in \boldsymbol L^{\alpha'} (\Omega)$. Hence, the definition of $a$ is sensible.
}
\item{
Since $\Omega$ has finite measure, $L^{\alpha}(\Omega) \subset L^1(\Omega)$. In particular, $\vert \nabla y \vert \in L^1(\Omega)$ which guarantees that the integral in $j$ is well-posed.
}
\end{enumerate}
\end{remark}

\begin{definition}
Let $f \in L^{\alpha'}(\Omega)$ be a given force density. We call a velocity field $y \in Y$ a \emph{weak solution to the Herschel-Bulkley problem} if
\begin{equation}
a(y,z - y) + j(z) - j(y) \geq \left( f , z - y \right)_{\alpha}, \quad \forall z \in Y. \label{eq-dl-vi}
\end{equation}
\end{definition}

It turns out that this variational inequality of the second kind is equivalent to the first order necessary optimality condition for the minimisation problem
\begin{equation}
\inf_{y \in Y} I(y) = c(y) + j(y) - \left( f , y \right)_{\alpha},\label{eq-dl-mp}
\end{equation}
cf.  \cite{Huilgol2005}, \cite[p~8]{Glowinski1981}.

Huilgol and You \cite{Huilgol2005} conclude that ``there are gaps as far as the existence and uniqueness of a solution to the pipe flows of the three fluids under consideration" \cite[p~141]{Huilgol2005} (the authors consider Bingham, Herschel-Bulkley and Casson fluids). In fact, it turns out that both existence and uniqueness can be established with well-known tools from convex optimisation:

\begin{theorem} \label{thm-dl-existence}
Problems \eqref{eq-dl-vi} and \eqref{eq-dl-mp} are equivalent and possess a unique solution $y^* \in Y$.
\begin{proof}
\emph{Existence.} A standard proof for the existence of a solution to \eqref{eq-dl-mp} can be found e.g. in Theorem 1.1 and Remark 1.2 of \cite[pp~7-8]{Lions1971}. It requires weak lower semi-continuity of the objective $I$ and
\begin{equation}
I(y) \to + \infty \quad \text{as } \Vert y \Vert_{W^{1,\alpha}(\Omega)} \to \infty \label{eq-dl-funcgrowth}
\end{equation}
to ensure boundedness of any minimising sequence.

The first assumption follows immediately since $I$ is both continuous and convex, cf. Theorem 3.3.3 in \cite[p~93]{Buttazzo2006}.

To verify the second assumption, we use the following generalisation of Poincar{\'e}'s inequality: there exists a constant $C>0$ independent of $y \in Y$ such that
\begin{equation}
\Vert y \Vert_{W^{1,\alpha}(\Omega)}^{\alpha} \leq C \left( \int\limits_\Omega \vert \nabla y \vert^{\alpha} \;\mathrm{d}x + \left( \int\limits_{\Gamma_\mathrm{D}} y \;\mathrm{d}s \right)^{\alpha} \right). \label{eq-dl-genpoinc}
\end{equation}
A proof for the special case $\alpha = 2$ is given in \cite[p~355]{Wloka1987}. An analysis shows that it remains valid for $1 < \alpha < \infty$.

Since $y=0$ almost everywhere on $\Gamma_\mathrm{D}$, the last term in \eqref{eq-dl-genpoinc} vanishes. Consequently, using \eqref{eq-dl-genpoinc} along with H{\"o}lder's inequality, we obtain
\begin{align*}
I(y) &= \frac{\kappa}{\alpha} \int\limits_\Omega \left\vert \nabla y \right\vert^{\alpha} \,\mathrm{d}x + \tau_0 \int\limits_\Omega \left\vert \nabla y \right\vert \,\mathrm{d}x - (f,y)_{\alpha}\\
&\geq C \Vert y \Vert_{W^{1,\alpha}(\Omega)}^{\alpha} + 0 - \Vert f \Vert_{\alpha'} \Vert y \Vert_{\alpha}\\
&\geq C \left( \Vert y \Vert_{W^{1,\alpha}(\Omega)}^{\alpha} - \Vert y \Vert_{W^{1,\alpha}(\Omega)} \right)\\
&\to + \infty \quad \text{as } \Vert y \Vert_{W^{1,\alpha}(\Omega)} \to \infty.
\end{align*}

\emph{Uniqueness.} As a consequence of the strict convexity of $I$, a solution of \eqref{eq-dl-mp} is unique.

\emph{Equivalence.} We refer to Theorem 1.6 and Remark 1.12 in \cite[pp~12-13]{Lions1971}.
\end{proof}
\end{theorem}

\subsection{Constrained formulation}

Augmented Lagrangian methods are a common choice for the numerical solution of nonsmooth optimisation problems. One first introduces a new variable $\boldsymbol q$ for the rate of strain $\nabla y$, where $\boldsymbol q$ and $\nabla y$ are linked in a constraint. This trick allows to decouple nonsmoothness and nonlinearity on the one hand from the linear velocity term on the other hand. The constraint is then relaxed and penalised in an augmented Lagrangian.

The unconstrained problem \eqref{eq-dl-mp} rewritten as artificially constrained problem becomes
\begin{subequations} \label{eq-alg2}
\begin{align}
\min_{(y,\boldsymbol q) \in Y\times\boldsymbol L^\alpha (\Omega)} \tilde{I}(y,\boldsymbol q) &= \frac{\kappa}{\alpha} \int\limits_\Omega \vert \boldsymbol q \vert^\alpha \;\mathrm{d}x + \tau_0 \int\limits_\Omega \vert \boldsymbol q \vert \;\mathrm{d}x - (f,y)_\alpha \label{eq-alg2-obj}\\
\text{subject to } & \boldsymbol q - \nabla y = 0. \label{eq-alg2-con}
\end{align}
\end{subequations}
We define the corresponding Lagrangian $L: Y\times \boldsymbol L^\alpha (\Omega) \times \boldsymbol L^{\alpha'} (\Omega)\to \mathbb{R}$ by
\begin{equation} \label{eq-lagr}
L(y,\boldsymbol q, \boldsymbol \tau) = \frac{\kappa}{\alpha} \int\limits_\Omega \vert \boldsymbol q \vert^\alpha \;\mathrm{d}x + \tau_0 \int\limits_\Omega \vert \boldsymbol q \vert \;\mathrm{d}x - (f,y)_\alpha - (\boldsymbol \tau, \boldsymbol q - \nabla y)_\alpha.
\end{equation}
We will now answer the question for the existence of Lagrange multipliers.

\begin{proposition}\label{thm-alg2-ex-lm}
Let $(y^*,\boldsymbol q^*)\in Y\times \boldsymbol L^\alpha(\Omega)$ be the solution of \eqref{eq-alg2}. There exists a Lagrange multiplier $\boldsymbol \tau^* \in \boldsymbol L^{\alpha'}(\Omega)$ such that the KKT conditions
\begin{subequations}
\begin{align}
D_y L(y^*,\boldsymbol q^*,\boldsymbol \tau^*) &= 0\\
\partial_{\boldsymbol q} L(y^*,\boldsymbol q^*,\boldsymbol \tau^*) &\ni 0\\
D_{\boldsymbol \tau} L(y^*,\boldsymbol q^*,\boldsymbol \tau^*) &= 0
\end{align}
\end{subequations}
hold, i.e.
\begin{subequations}
\begin{align}
(\boldsymbol \tau^*, \nabla z)_\alpha &= (f,z)_\alpha \quad \forall z \in Y\\
\left(\kappa  \vert \boldsymbol q^* \vert^{\alpha-2} \boldsymbol q^* - \boldsymbol \tau^*,\boldsymbol r - \boldsymbol q^*\right)_\alpha + \tau_0 \int\limits_\Omega \vert \boldsymbol r \vert \;\mathrm{d}x - \tau_0 \int\limits_\Omega \vert \boldsymbol q^* \vert \;\mathrm{d}x &\geq 0 \quad \forall \boldsymbol r \in \boldsymbol L^\alpha (\Omega)\\
\boldsymbol q^* &= \nabla y^* \quad \text{a.e. in } \Omega.
\end{align}
\end{subequations}
Here, $D_y$ denotes the Fr\'echet derivative with respect to $y$ and $\partial_{\boldsymbol q}$ the subdifferential with respect to $\boldsymbol q$.
\begin{proof}
The assertion is an immediate consequence of the fact that Slater's constraint qualification (SCQ) is trivially satisfied here (cf \cite[Thm.~3.34]{Ruszczynski2006}).
\end{proof}
\end{proposition}
Physically, a Lagrange multiplier $\boldsymbol \tau$ associated with the constraint \eqref{eq-alg2-con} can be interpreted as an admissible stress \cite{Huilgol2005}. Therefore, the (Lagrangian) dual problem to \eqref{eq-alg2} can be seen as an optimisation problem in terms of the stress $\boldsymbol \tau$ instead of the veloctiy $y$. We will now follow this dual approach.

\section{Dual formulation}
\label{sec-dual}

\subsection{The Lagrangian dual}

Lagrangian duality expresses the fact that a constrained optimisation problem can be represented in two distinct, but often equivalent ways: the primal problem in terms of the primal variables or the dual problem in terms of the dual variables (Lagrange multipliers).

From the Lagrangian $L$ in \eqref{eq-lagr}, one obtains a primal function $L_\mathrm{p}$ by maximising in the dual variable
\begin{align*}
L_\mathrm{p}(y,\boldsymbol q) &:= \sup_{\boldsymbol \tau \in \boldsymbol L^{\alpha'}(\Omega)} L(y,\boldsymbol q, \boldsymbol \tau)\\
\intertext{and a dual function $L_\mathrm{d}$ by minimising in the primal variables}
L_\mathrm{d}(\boldsymbol \tau) &:= \inf_{(y,\boldsymbol q)\in Y\times \boldsymbol L^\alpha(\Omega)} L(y,\boldsymbol q, \boldsymbol \tau).
\end{align*}
The primal and dual functions are permitted to assume values on the extended real line, $\mathbb{R} \cup \lbrace +\infty \rbrace$ and $\mathbb{R} \cup \lbrace -\infty \rbrace$, respectively. Then the primal problem is given by
\begin{equation*}
\min_{(y,\boldsymbol q)\in Y\times \boldsymbol L^\alpha(\Omega)} L_\mathrm{p}(y,\boldsymbol q)
\end{equation*}
and the dual problem by
\begin{equation*}
\max_{\boldsymbol \tau \in \boldsymbol L^{\alpha'}(\Omega)} L_\mathrm{d}(\boldsymbol \tau).
\end{equation*}
Since we have just seen that SCQ holds, we actually have strong duality
\begin{equation*}
\min_{(y,\boldsymbol q)\in Y\times \boldsymbol L^\alpha(\Omega)} L_\mathrm{p}(y,\boldsymbol q) = \max_{\boldsymbol \tau \in \boldsymbol L^{\alpha'}(\Omega)} L_\mathrm{d}(\boldsymbol \tau),
\end{equation*}
i.e. there is no duality gap. Furthermore, with solutions $(y^*,\boldsymbol q^*)$ and $\boldsymbol \tau^*$ of the primal and dual problem, respectively, $(y^*,\boldsymbol q^*,\boldsymbol \tau^*)$ is a saddle point of the Lagrangian (primal-dual solution), see e.g. \cite[Thms.~4.7, 4.9]{Ruszczynski2006}:
\begin{equation*}
L(y^*,\boldsymbol q^*, \boldsymbol \tau) \leq L(y^*,\boldsymbol q^*, \boldsymbol \tau^*) \leq L(y,\boldsymbol q, \boldsymbol \tau^*) \quad \forall (y,\boldsymbol q,\boldsymbol \tau)\in Y\times \boldsymbol L^\alpha(\Omega) \times \boldsymbol L^{\alpha'}(\Omega).
\end{equation*}

Calculating the primal function from the Lagrangian $L$ in \eqref{eq-lagr} yields
\begin{equation*}
L_\mathrm{p}(y,\boldsymbol q) = 
\begin{cases}
\frac{\kappa}{\alpha} \int\limits_\Omega \vert \boldsymbol q \vert^\alpha \;\mathrm{d}x + \tau_0 \int\limits_\Omega \vert \boldsymbol q \vert \;\mathrm{d}x - (f,y)_\alpha & \text{if } \boldsymbol q = \nabla y\\
+\infty & \text{otherwise}
\end{cases}
\end{equation*}
Clearly, the primal problem $\min L_\mathrm{p} (y,\boldsymbol q)$ is equivalent to \eqref{eq-alg2}.

Let us now turn to the dual problem. We assume $\boldsymbol \tau \in \boldsymbol L^{\alpha'}(\Omega)$ is given and we will now calculate the corresponding value of the dual function $L_\mathrm{d}(\boldsymbol \tau)$.

Initially we observe, since $L$ depends linearly on $y$, that $L_\mathrm{d}(\boldsymbol \tau) = -\infty$ if $(\boldsymbol \tau, \nabla \cdot)_\alpha \neq (f,\cdot)_\alpha$ in $Y^*$. This is a weak formulation of $-\mathrm{div}\, \boldsymbol \tau \neq f$, stating that $\boldsymbol \tau$ violates conservation of momentum \eqref{eq-hb-c-com}.

Otherwise, if $(\boldsymbol \tau, \nabla \cdot)_\alpha = (f,\cdot)_\alpha$ holds as equation in $Y^*$, we have to find $\boldsymbol q^* \in \boldsymbol L^\alpha (\Omega)$ which minimises the remaining convex functional
\begin{equation} \label{eq-alg2-trs-q-tau}
\boldsymbol q^* = \argmin_{\boldsymbol q \in \boldsymbol L^\alpha (\Omega)} \frac{\kappa}{\alpha} \int\limits_\Omega \vert \boldsymbol q \vert^\alpha \;\mathrm{d}x + \tau_0 \int\limits_\Omega \vert \boldsymbol q \vert \;\mathrm{d}x - (\boldsymbol \tau , \boldsymbol q)_\alpha .
\end{equation}

\begin{lemma}
The unique solution of \eqref{eq-alg2-trs-q-tau} is
\begin{equation}
\boldsymbol q^* = \frac{1}{\kappa^{1/(\alpha-1)}} (\vert \boldsymbol \tau \vert - \tau_0)_+^{1/(\alpha-1)} \frac{\boldsymbol \tau}{\vert \boldsymbol \tau \vert},
\end{equation}
where $(\cdot)_+ := \max(0,\cdot)$ represents the hard-thresholding operator for truncation below zero.
\begin{proof}
\emph{Feasibility.} First of all, we show that the above $\boldsymbol q^*$ is admissible, i.e. belongs to $\boldsymbol L^\alpha (\Omega)$. Since $\tau \in \boldsymbol L^{\alpha'} (\Omega)$, we have $\vert \boldsymbol \tau \vert - \tau_0 \in L^{\alpha'} (\Omega)$ and thus
\begin{align*}
\left\Vert \left( \vert \boldsymbol \tau \vert - \tau_0 \right)_+^{\frac{1}{\alpha - 1}} \frac{\boldsymbol \tau}{\vert \boldsymbol \tau \vert} \right\Vert_{\alpha}^{\alpha} &= \int\limits_\Omega \left\vert \left( \vert \boldsymbol \tau \vert_2 - \tau_0 \right)_+^{\frac{1}{\alpha - 1}} \frac{\boldsymbol \tau}{\vert \boldsymbol \tau \vert_2} \right\vert_{\alpha}^{\alpha} \,\mathrm{d}x\\
&\leq C \int\limits_\Omega \left\vert \left( \vert \boldsymbol \tau \vert_2 - \tau_0 \right)_+^{\frac{1}{\alpha - 1}} \frac{\boldsymbol \tau}{\vert \boldsymbol \tau \vert_2} \right\vert_2^{\alpha} \,\mathrm{d}x\\
&= C \int\limits_\Omega \left( \vert \boldsymbol \tau \vert_2 - \tau_0 \right)_+^{\frac{\alpha}{\alpha - 1}} \frac{\vert \boldsymbol \tau \vert_2^{\alpha}}{\vert \boldsymbol \tau \vert_2^{\alpha}} \,\mathrm{d}x\\
&= C \int\limits_\Omega \left( \vert \boldsymbol \tau \vert - \tau_0 \right)_+^{\alpha'} \,\mathrm{d}x\\
&\leq C \int\limits_\Omega \left\vert \vert \boldsymbol \tau \vert - \tau_0 \right\vert^{\alpha'} \,\mathrm{d}x < \infty.
\end{align*}
Consequently, $\left( \vert \boldsymbol \tau \vert - \tau_0 \right)_+^{\frac{1}{\alpha - 1}} \frac{\boldsymbol \tau}{\vert \boldsymbol \tau \vert} \in \boldsymbol L^{\alpha} (\Omega)$.

\emph{Optimality.} In analogy to Proposition \ref{thm-alg2-ex-lm} we conclude that the minimising argument satisfies the first-order optimality condition
\begin{equation} \label{eq-alg2-q-vi}
\left(\kappa  \vert \boldsymbol q^* \vert^{\alpha-2} \boldsymbol q^* - \boldsymbol \tau,\boldsymbol r - \boldsymbol q^*\right)_\alpha + \tau_0 \int\limits_\Omega \vert \boldsymbol r \vert \;\mathrm{d}x - \tau_0 \int\limits_\Omega \vert \boldsymbol q^* \vert \;\mathrm{d}x \geq 0 \quad \forall \boldsymbol r \in \boldsymbol L^\alpha (\Omega).
\end{equation}
Due to the strict convexity of the problem, this condition provides a characterisation for the unique solution $\boldsymbol q^*$.

We set $\boldsymbol s := \boldsymbol r - \boldsymbol q^*$. Using the above $\boldsymbol q^*$ we obtain for the first term
\begin{align*}
\left(\kappa  \vert \boldsymbol q^* \vert^{\alpha-2} \boldsymbol q^* - \boldsymbol \tau,\boldsymbol s\right)_\alpha &=  \left( ( \vert \boldsymbol \tau \vert - \tau_0 )_+ \frac{\boldsymbol \tau}{\vert \boldsymbol \tau \vert} - \boldsymbol \tau, \boldsymbol s \right)_\alpha\\
&= \int\limits_{\lbrace\vert\boldsymbol\tau\vert>\tau_0\rbrace} ( \vert \boldsymbol \tau \vert - \tau_0 ) \frac{\boldsymbol \tau}{\vert \boldsymbol \tau \vert} \cdot \boldsymbol s \, \mathrm{d}x - (\boldsymbol \tau, \boldsymbol s)_\alpha\\
&= - \int\limits_{\lbrace\vert\boldsymbol\tau\vert\leq\tau_0\rbrace} \boldsymbol \tau \cdot \boldsymbol s \, \mathrm{d}x - \tau_0\int\limits_{\lbrace\vert\boldsymbol\tau\vert>\tau_0\rbrace} \frac{\boldsymbol \tau}{\vert \boldsymbol \tau \vert} \cdot \boldsymbol s \, \mathrm{d}x.
\end{align*}
These two integrals can be estimated from below: from the Cauchy-Schwarz inequality we derive
\begin{equation}
\int\limits_{\lbrace\vert\boldsymbol\tau\vert\leq\tau_0\rbrace} \boldsymbol \tau \cdot \boldsymbol s \, \mathrm{d}x \leq  \int\limits_{\lbrace\vert\boldsymbol\tau\vert\leq\tau_0\rbrace} \vert \boldsymbol \tau \vert \vert \boldsymbol s \vert \, \mathrm{d}x \leq  \tau_0 \int\limits_{\lbrace\vert\boldsymbol\tau\vert\leq\tau_0\rbrace} \vert \boldsymbol s \vert \, \mathrm{d}x = \tau_0 \int\limits_{\lbrace\vert\boldsymbol\tau\vert\leq\tau_0\rbrace} \vert \boldsymbol s + \boldsymbol q^* \vert \, \mathrm{d}x,
\end{equation}
since $\boldsymbol q^* = 0$ almost everywhere on $\lbrace\vert\boldsymbol\tau\vert\leq\tau_0\rbrace$. Similarly we conclude
\begin{align*}
\tau_0\int\limits_{\lbrace\vert\boldsymbol\tau\vert>\tau_0\rbrace} \frac{\boldsymbol \tau}{\vert \boldsymbol \tau \vert} \cdot \boldsymbol s \, \mathrm{d}x
&= \tau_0\int\limits_{\lbrace\vert\boldsymbol\tau\vert>\tau_0\rbrace} \frac{\boldsymbol \tau}{\vert \boldsymbol \tau \vert} \cdot \boldsymbol s \, \mathrm{d}x + \tau_0 \int\limits_{\lbrace\vert\boldsymbol\tau\vert>\tau_0\rbrace} \vert \boldsymbol q^* \vert \,\mathrm{d}x - \tau_0 \int\limits_\Omega \vert \boldsymbol q^* \vert \, \mathrm{d}x\\
&= \tau_0\int\limits_{\lbrace\vert\boldsymbol\tau\vert>\tau_0\rbrace} \frac{\boldsymbol \tau}{\vert \boldsymbol \tau \vert} \cdot \left( \boldsymbol s + \vert \boldsymbol q^* \vert \frac{\boldsymbol \tau}{\vert \boldsymbol \tau \vert} \right) \,\mathrm{d}x - \tau_0 \int\limits_\Omega \vert \boldsymbol q^* \vert \, \mathrm{d}x\\
&\leq \tau_0\int\limits_{\lbrace\vert\boldsymbol\tau\vert>\tau_0\rbrace} \left\vert \boldsymbol s + \vert \boldsymbol q^* \vert \frac{\boldsymbol \tau}{\vert \boldsymbol \tau \vert} \right\vert \,\mathrm{d}x - \tau_0 \int\limits_\Omega \vert \boldsymbol q^* \vert \, \mathrm{d}x\\
&= \tau_0\int\limits_{\lbrace\vert\boldsymbol\tau\vert>\tau_0\rbrace} \vert \boldsymbol s + \boldsymbol q^* \vert \,\mathrm{d}x - \tau_0 \int\limits_\Omega \vert \boldsymbol q^* \vert \, \mathrm{d}x,
\end{align*}
where we used $\boldsymbol q^* / \vert \boldsymbol q^* \vert = \boldsymbol \tau / \vert \boldsymbol \tau \vert$ almost everywhere on $\lbrace \vert\boldsymbol\tau\vert>\tau_0\rbrace$, which proves the assertion.
\end{proof}
\end{lemma}

By substituting this solution $\boldsymbol q^*$ in the Lagrangian \eqref{eq-lagr} and simplifying, we obtain the dual function
\begin{equation*}
L_\mathrm{d}(\boldsymbol \tau) = 
\begin{cases}
-\frac{1}{\alpha'\kappa^{1/(\alpha-1)}} \int\limits_\Omega (\vert \boldsymbol \tau \vert- \tau_0)_+^{\alpha'} \;\mathrm{d}x & \text{if } (\boldsymbol \tau, \nabla \cdot)_\alpha = (f,\cdot)_\alpha \text{ in } Y^*\\
-\infty & \text{otherwise}
\end{cases}
\end{equation*}
Instead of the dual problem $\max L_\mathrm{d} (\boldsymbol \tau)$, we will from now on consider the equivalent problem
\begin{subequations} \label{eq-tt-mp}
\begin{align}
\min_{\boldsymbol \tau \in \boldsymbol L^{\alpha'} (\Omega)} J(\boldsymbol \tau) &= \frac{1}{\alpha'\kappa^{1/(\alpha-1)}} \int\limits_\Omega (\vert \boldsymbol \tau \vert- \tau_0)_+^{\alpha'} \;\mathrm{d}x \label{eq-tt-mp-obj}\\
\text{subject to } & (\boldsymbol \tau, \nabla \cdot)_\alpha = (f,\cdot)_\alpha \text{ in } Y^* \label{eq-tt-mp-con}
\end{align}
\end{subequations}
and, with a slight abuse of terminology, refer to \eqref{eq-tt-mp} as the dual problem.

The numerical method we propose will be based on this characterisation of Herschel-Bulkey duct flow. We remark that \eqref{eq-tt-mp} is convex, but generally not strictly convex. Hence, while a solution $\boldsymbol \tau^* \in \boldsymbol L^{\alpha'}(\Omega)$ is always guaranteed to exist, it needs not be unique.

\subsection{Optimality conditions}

We already saw, since SCQ holds for the primal problem, that the primal and dual problem possess solutions and that strong duality holds, i.e. both problems are equivalent. We could also arrive at this conclusion by verifying, say, the linear independence constraint qualification (LICQ) for the dual problem, which for any dual solution $\boldsymbol \tau^* \in \boldsymbol L^{\alpha'}(\Omega)$ implies the existence of a (unique) Lagrange multiplier $y^*\in Y$.

Due to the convexity of problem \eqref{eq-tt-mp}, the corresponding KKT conditions
\begin{subequations}
\label{eq-tt-os}
\begin{align}
\frac{1}{\kappa^{1/(\alpha-1)}} \left( \boldsymbol \sigma, \left( \vert \boldsymbol \tau \vert - \tau_0 \right)_+^{\frac{1}{\alpha - 1}} \frac{\boldsymbol \tau}{\vert \boldsymbol \tau \vert} \right)_\alpha - \left( \boldsymbol \sigma , \nabla y \right)_\alpha = 0 \quad &\forall \boldsymbol \sigma \in \boldsymbol L^{\alpha'}(\Omega) \label{eq-tt-os-vp}\\
\left( \boldsymbol \tau, \nabla z \right)_\alpha = \left( f , z \right)_{\alpha} \quad &\forall z \in Y. \label{eq-tt-os-com}
\end{align}
\end{subequations}
are necessary and sufficient for optimality. This justifies the following definition:

\begin{definition} \label{def-mixed}
We call a pair $(\boldsymbol \tau,y) \in \boldsymbol L^{\alpha'} (\Omega) \times Y$ a \emph{weak solution to the mixed Herschel-Bulkley problem}, if they satisfy \eqref{eq-tt-os}.
\end{definition}

One can relate \eqref{eq-tt-os} back to the strong formulation \eqref{eq-hb-c} we started with: looking at \eqref{eq-tt-os-vp}, a corresponding strong form reads
\begin{equation*}
\frac{1}{\kappa^{1/(\alpha-1)}} \left( \vert \boldsymbol \tau \vert - \tau_0 \right)_+^{\frac{1}{\alpha - 1}} \frac{\boldsymbol \tau}{\vert \boldsymbol \tau \vert} - \nabla y = 0.
\end{equation*}
In plain words, the rate of strain can be recovered from the stress by a soft-thresholding operation, combined with a power law if $\alpha \neq 2$. By rearranging for the stress $\boldsymbol \tau$, we see that this is equivalent to the viscoplastic constitutive relations of the Herschel-Bulkley model in \eqref{eq-hb-c-visc}-\eqref{eq-hb-c-plas}, cf \cite{Saramito2007,Saramito2009}:
\begin{equation*}
\nabla y = \frac{1}{\kappa^{1/(\alpha-1)}} \left( \vert \boldsymbol \tau \vert - \tau_0 \right)_+^{\frac{1}{\alpha - 1}} \frac{\boldsymbol \tau}{\vert \boldsymbol \tau \vert} \Leftrightarrow 
\begin{cases}
\boldsymbol\tau = \kappa \left\vert \nabla y \right\vert^{\alpha - 2} \nabla y + \tau_0 \frac{\nabla y}{\left\vert \nabla y \right\vert} & \mbox{if } \nabla y \neq 0\\
\left\vert {\boldsymbol\tau} \right\vert \leq \tau_0 & \mbox{if } \nabla y = 0.
\end{cases}
\end{equation*}
We already identified \eqref{eq-tt-os-com} as a weak formulation of conservation of momentum \eqref{eq-hb-c-com}. The boundary conditions \eqref{eq-hb-c-dbc}-\eqref{eq-hb-c-nbc} are incorporated through an appropriate choice of function spaces for the solution and the test functions.

\subsection{Regularity of the objective}
\label{sec-regularity}

Numerical optimisation algorithms typically benefit in terms of their speed of convergence when higher derivatives of the objective and any constraints exist and are used. Due to the nonsmoothness in $j$, the conventional functional $I$ in \eqref{eq-dl-mp} is not Fr{\'e}chet differentiable whenever a vanishing velocity gradient occurs. In contrast, the stress functional $J$ in \eqref{eq-tt-mp} is generally continuously Fr{\'e}chet differentiable in $\boldsymbol \tau$. Depending on the value of $\alpha$, additional regularity may be achieved.

For shear-thickening Herschel-Bulkley fluids, $\alpha > 2$, no higher Fr{\'e}chet derivatives of the objective $J$ are available. This is due to the fact that the optimality system \eqref{eq-tt-os} is not differentiable at the interface between viscous and plastic regions:

Let $\boldsymbol \tau_* \in \boldsymbol L^{\alpha}(\Omega)$ with $\vert \boldsymbol \tau_* \vert = \tau_0$ be the limit of a sequence $(\boldsymbol \tau^n)_{n\in \mathbb{N}}$ with $\vert \boldsymbol \tau^n \vert < \tau_0$ for all $n\in\mathbb{N}$ and $\vert \boldsymbol \tau^n \vert \nearrow \tau_0$ as $n\to\infty$. Let $(\boldsymbol \sigma^n)_{n\in \mathbb{N}}$ be another sequence with $\lim_{n\to\infty} \boldsymbol \sigma_n = \boldsymbol \tau_*$, $\vert \boldsymbol \sigma^n \vert > \tau_0$ for all $n\in\mathbb{N}$ and $\vert \boldsymbol \sigma^n \vert \searrow \tau_0$ as $n\to\infty$. $J$ is twice Fr{\'e}chet differentiable at every $\boldsymbol \tau^n$ and every $\boldsymbol \sigma^n$. From \eqref{eq-tt-mp-obj} we observe, however, that
\begin{align*}
\lim_{n\to\infty} \Vert J''(\boldsymbol \tau^n) \Vert &= 0
\intertext{while}
\lim_{n\to\infty} \Vert J''(\boldsymbol \sigma^n) \Vert &= \infty.
\end{align*}
This sudden jump from an infinite to a vanishing viscosity is not a flaw of the analytical framework chosen here, but rather an unphysical feature of the model.

For shear-thinning Herschel-Bulkley fluids, $1 < \alpha < 2$, one can conclude from \eqref{eq-tt-mp-obj} that $J$ is at least twice continuously Fr{\'e}chet differentiable. 

For Bingham fluids, $\alpha = 2$, $J$ is once, but not twice continuously Fr{\'e}chet differentiable. However, $J'$ still satisfies a slightly weaker notion of differentiability that is often applied in the context of semismooth Newton methods.
 
\begin{definition}
Let $K,L$ be Banach spaces, $D\subset K$ an open subset.

A function $F: K \supset D \to L$ is called \emph{Newton differentiable} (or \emph{slantly differentiable}) in the open subset $U \subset D$, if there exists a family of mappings $G: U \to \mathcal{L}(K,L)$ such that
\begin{equation}
\lim_{h\to 0} \frac{\left\Vert F(x+h) - F(x) - G(x+h)h \right\Vert_L}{\left\Vert h \right\Vert_K} = 0, \quad \forall x\in U.
\end{equation}
\end{definition}

For a proof that $J'$ is indeed Newton differentiable, we refer to Lemma 3.1 in \cite{Hintermueller2002}.

This result has far-reaching consequences: Sequential Quadratic Programming (SQP) methods attempt to solve the system of first order necessary optimality conditions with methods of Newton type. Since, independent of $\alpha$, the functional $I$ does not possess a continuous Fr{\'e}chet derivative where $\nabla y = 0$, in particular no second Fr{\'e}chet or Newton derivative, SQP methods are not applicable to this problem. With the alternative approach suggested here that minimises $J$, Newton differentiability (or even higher regularity) of the optimality system is warranted for the most important applications with shear-thinning and Bingham fluids, $1 < \alpha \leq 2$. Therefore, the new approach to formulating viscoplastic flow allows us to use powerful numerical algorithms from the SQP framework.

A similar duality approach was chosen by de los Reyes and Gonz\'alez Andrade \cite{Reyes2009,Reyes2010} for Bingham fluids. However, the authors used a different splitting in the primal objective, which results in a more complicated dual problem. They then introduce a Tikhonov regularisation into the dual objective and use a semismooth Newton method to solve the resulting optimality system. We could certainly follow an analogous procedure here and solve a penalised dual problem with a semismooth Newton method for Bingham flow, and a smooth Newton method for shear-thinning Herschel-Bulkley flow. However, the very simple structure of the dual \eqref{eq-tt-mp} allows us to even solve the unregularised problem. In our numerical method, this will require the application of trust regions, resulting in exactly viscoplastic solutions.

\section{Discretisation}
\label{sec-discretisation}

For the numerical solution of \eqref{eq-tt-mp}, we follow a first-discretise-then-optimise approach. That is, we first discretise both the objective functional and the constraint in order to obtain an optimisation problem that is posed in finite dimensions. We are then in the position to apply well-established methods to solve this problem numerically.

Let $\Omega^h$ be a polygonal domain that approximates $\Omega$, $\emptyset \neq \Gamma_\mathrm{D}^h \subset \partial \Omega^h$ the set of all edge segments that correspond to $\Gamma_\mathrm{D}$, analogously $\Gamma_\mathrm{N}^h = \partial \Omega^h \setminus \Gamma_\mathrm{D}^h$ and $\mathcal{T}^{h}$ a regular triangulation of the closure $\bar{\Omega}^h$. We discretise with P1-P0 finite elements, i.e. discretise the stress with piecewise constant functions, the velocity with piecewise linear functions. Denoting with $P_k$ the space of polynomials in $\mathbb{R}^2$ with degree $k$ or less, we therefore introduce the following finite-dimensional spaces:
\begin{align*}
\boldsymbol S^h &:= \left\lbrace \boldsymbol \tau^h = (\tau_1^h,\tau_2^h) \in \boldsymbol L^{\alpha'}(\Omega^h) : \tau_j^h\vert_T \in P_0, \; \forall T \in \mathcal{T}^h, j \in \lbrace 1,2 \rbrace \right\rbrace \subset \boldsymbol L^{\alpha'}(\Omega^h)\\
V^h &:= \left\lbrace y^h \in C(\Omega^h) : y^h\vert_T \in P_1, \; \forall T \in \mathcal{T}^h \right\rbrace \subset W^{1,\alpha}(\Omega^h).
\intertext{Furthermore, we define the subspace of discrete velocity functions}
Y^h &:= \left\lbrace y^h \in V^h : y^h = 0 \text{ on } \Gamma_\mathrm{D}^h \right\rbrace \subset V^h.
\end{align*}

A basis of the space $\boldsymbol S^h$ is given by
\begin{equation}
\boldsymbol \chi_i :=
\begin{cases}
\chi_k \vec{e}_1 & \text{if } i = 2k-1\\
\chi_k \vec{e}_2 & \text{if } i = 2k\\
\end{cases}
\quad i = 1,...,2n_\mathrm{T}.
\end{equation}
Here, $\chi_k$ is the characteristic function of the triangle $T_k\in\mathcal{T}^h$, $n_\mathrm{T}$ is the number of triangles in $\mathcal{T}^h$ and $\vec{e}_1, \vec{e}_2$ are the canonical unit vectors of $\mathbb{R}^2$.

In the following we will assume that the $n_\mathrm{N}$ triangle nodes $(x^j)_{j=1}^{n_\mathrm{N}}$ are sorted such that the first $n_\mathrm{I}$ nodes, $n_\mathrm{I} < n_\mathrm{N}$, lie in the interior or on the Neumann boundary of the domain, $(x^j)_{j=1}^{n_\mathrm{I}} \subset \Omega^h \cup \Gamma_\mathrm{N}^h$, the remaining $n_\mathrm{D} = n_\mathrm{N}-n_\mathrm{I}$ nodes on the Dirichlet boundary, $(x^j)_{j=n_\mathrm{I}+1}^{n_\mathrm{N}} \subset \Gamma_\mathrm{D}^h$. By $\phi_j, j = 1,...,n_\mathrm{N}$ we denote the hat function of the node $x^j$. Then $(\phi_j)_{j=1}^{n_\mathrm{N}}$ form a basis for $V^h$ and $(\phi_j)_{j=1}^{n_\mathrm{I}}$ span the subspace $Y^h$.

Using these spaces, a discrete counterpart of \eqref{eq-tt-mp} reads as follows: find a solution $\boldsymbol \tau^h \in \boldsymbol S^h$ for the minimisation problem
\begin{subequations} \label{eq-tt-hb-dsc-mp}
\begin{align}
\inf_{\boldsymbol \tau^h \in \boldsymbol S^h} J(\boldsymbol\tau^h) &= \frac{1}{\alpha'\kappa^{1/(\alpha-1)}} \int\limits_\Omega \left( \vert \boldsymbol \tau^h \vert - \tau_0 \right)_+^{\alpha'} \;\mathrm{d}x \label{eq-tt-hb-dsc-mp-obj}\\
\text{subject to } &\left( \boldsymbol \tau^h, \nabla \cdot \right)_\alpha = \left( f , \cdot \right)_{\alpha} \quad \text{in } (Y^h)^* \label{eq-tt-hb-dsc-mp-com}
\end{align}
\end{subequations}

An equivalent representation of \eqref{eq-tt-hb-dsc-mp} can be obtained through a standard procedure, which yields an optimisation problem posed in $\mathbb{R}^{2n_\mathrm{T}}$:
\begin{subequations} \label{eq-tt-hb-fds-mp}
\begin{align}
\inf J^h(\vec{\boldsymbol\tau}) &= \frac{1}{\alpha'\kappa^{1/(\alpha-1)}} \sum_{k=1}^{n_\mathrm{T}} \vert T_k \vert
\left(
\vert \vec{\boldsymbol \tau}_k \vert - \tau_0
\right)_+^{\alpha'} \label{eq-tt-hb-fds-mp-obj}\\
\text{subject to } &D^h \vec{\boldsymbol \tau} = f^h \label{eq-tt-hb-fds-mp-com}
\end{align}
\end{subequations}

$D^h \in \mathbb{R}^{n_\mathrm{I}\times 2n_\mathrm{T}}$, $f^h \in \mathbb{R}^{n_\mathrm{I}}$ and $\vec{\boldsymbol \tau}\in\mathbb{R}^{2n_\mathrm{T}}$ are the representations of $(\cdot,\nabla \cdot)_\alpha$, $(f,\cdot)_\alpha$ or $\boldsymbol \tau^h$ in terms of $(\boldsymbol \chi_i)_{i=1}^{2n_\mathrm{T}}$ or $(\phi_j)_{j=1}^{n_\mathrm{I}}$, respectively. We use $\vec{\boldsymbol \tau}_k$ as an abbreviation for $(\vec{\tau}_{2k-1},\vec{\tau}_{2k})^\top\in\mathbb{R}^2$.

\begin{proposition}
Problem \eqref{eq-tt-hb-fds-mp} has a solution $\vec{\boldsymbol \tau}^* \in \mathbb{R}^{2n_\mathrm{T}}$.
\begin{proof}
First, the admissible set for potential minimisers $\vec{\boldsymbol \tau}\in\mathbb{R}^{2n_\mathrm{T}}$ as defined through the constraint $D^h \vec{\boldsymbol \tau} = f^h$ is nonempty: since, due to the inf-sup stable definition of the finite element spaces, $D^h$ has full rank $n_\mathrm{I}$ \cite{Huilgol2005} and $D^h(D^h)^\top$ is invertible. One particular feasible point is given by
\begin{equation*}
\vec{\boldsymbol \tau} := (D^h)^\top (D^h(D^h)^\top)^{-1} f^h.
\end{equation*}

Furthermore, $J^h$ is convex, continuous and $J^h(\vec{\boldsymbol \tau}) \to +\infty$ as $\vert \vec{\boldsymbol \tau} \vert \to \infty$. In analogy to Theorem \ref{thm-dl-existence}, we conclude that there exists a minimiser $\vec{\boldsymbol \tau}^* \in \mathbb{R}^{2n_\mathrm{T}}$.
\end{proof}
\end{proposition}

\begin{proposition}
If $\vec{\boldsymbol \tau}^* \in \mathbb{R}^{2n_\mathrm{T}}$ is a solution to \eqref{eq-tt-hb-fds-mp}, then there exists a Lagrange multiplier $\vec{y}^* \in \mathbb{R}^{n_\mathrm{I}}$ such that $(\vec{\boldsymbol \tau}^*, \vec{y}^*)$ satisfies the KKT conditions
\begin{subequations} \label{eq-tt-kkt-fds}
\begin{align}
\nabla J^h (\vec{\boldsymbol \tau}^*) - (D^h)^\top \vec{y}^* &= 0\label{eq-tt-kkt1-fds}\\
D^h \vec{\boldsymbol \tau}^* = f^h,\label{eq-tt-kkt2-fds}
\end{align}
\end{subequations}
where $\nabla J^h (\vec{\boldsymbol \tau}^*) = (\nabla_{\vec{\boldsymbol \tau}_k} J^h (\vec{\boldsymbol \tau}^*))_{k=1}^{n_\mathrm{T}}$ is given by
\begin{equation} \label{eq-tt-grad}
\nabla_{\vec{\boldsymbol \tau}_k} J^h (\vec{\boldsymbol \tau}) = \frac{\vert T_k \vert}{\kappa^{1/(\alpha-1)}} \left(
\vert \vec{\boldsymbol \tau}_k \vert - \tau_0 \right)_+^{\frac{1}{\alpha-1}} \frac{\vec{\boldsymbol \tau}_k}{\vert \vec{\boldsymbol \tau}_k \vert}.
\end{equation}
\begin{proof}
This result follows immediately from the inf-sup stability of the chosen finite element spaces, i.e. $D^h$ is onto. Hence, LICQ holds for any $\vec{\boldsymbol \tau} \in \mathbb{R}^{n_\mathrm{T}}$.
\end{proof}
\end{proposition}

\section{Trust-region SQP algorithm}
\label{sec-algorithm}

We will now present the trust-region SQP algorithm that we propose to solve problem \eqref{eq-tt-hb-fds-mp}. It is based on the Byrd-Omojokun method \cite{Byrd1987,Omojokun1990} using an implementation similar to \cite{Lalee1998}.

In order to guarantee sufficient smoothness of the objective, we assume from now on $1 < \alpha \leq 2$. Then the Hessian $\nabla^2 J^h(\vec{\boldsymbol \tau})$ is locally bounded. From \eqref{eq-tt-grad} we obtain the following explicit representation:
\begin{equation} \label{eq-tt-hess}
\nabla_{\vec{\boldsymbol \tau}_k}^2 J^h (\vec{\boldsymbol \tau}) = \frac{\vert T_k \vert}{\kappa^{\frac{1}{\alpha-1}}(\alpha-1)} \frac{
\left(\vert\vec{\boldsymbol \tau}_k\vert - \tau_0\right)_+^{\frac{1}{\alpha-1}-1}
}{
\vert \vec{\boldsymbol \tau}_k \vert^3
} H(\vec{\boldsymbol \tau}_k),
\end{equation}
where the symmetric $2\times 2$ matrix $H$ has the components
\begin{align*}
h_{11} &= (\alpha-1) \vec{\tau}_{2k}^2 \left(\vert\vec{\boldsymbol \tau}_k\vert - \tau_0\right)_+ + \vec{\tau}_{2k-1}^2 \vert \vec{\boldsymbol \tau}_k \vert\\
h_{12} &= -\vec{\tau}_{2k-1}\vec{\tau}_{2k}\left((\alpha-1) \left(\vert\vec{\boldsymbol \tau}_k\vert - \tau_0\right)_+ - \vert \vec{\boldsymbol \tau}_k \vert \right)\\
h_{22} &= (\alpha-1) \vec{\tau}_{2k-1}^2 \left(\vert\vec{\boldsymbol \tau}_k\vert - \tau_0\right)_+ + \vec{\tau}_{2k}^2 \vert \vec{\boldsymbol \tau}_k \vert.
\end{align*}

\subsection{Sequential Quadratic Programming}

The rationale behind Sequential Quadratic Programming (SQP) is to approximate the nonlinear minimisation problem \eqref{eq-tt-hb-fds-mp} with a sequence of linear-quadratic problems, i.e. optimisation problems with a quadratic objective and linear constraints. For problems with equality constraints only, like here, a basic SQP method is equivalent to applying Newton's method to the first order optimality conditions \eqref{eq-tt-kkt-fds}. Given an iterate $(\vec{\boldsymbol{\tau}}^k,\vec{y}^k)$ ($k = 0, 1, 2, ...$), Newton's method attempts to find $(\vec{\boldsymbol{\delta\tau}},\vec{\delta y})$ such that
\begin{subequations} \label{eq-trs-newton}
\begin{align}
\nabla^2 J^k \vec{\boldsymbol{\delta\tau}} - D^\top \vec{\delta y} &= -(\nabla J^k - D^\top \vec{y}^k)\\
D \vec{\boldsymbol{\delta\tau}} &= -(D \vec{\boldsymbol{\tau}}^k - f).
\end{align}
\end{subequations}
To simplify the notation, we have suppressed the superscript $h$ and we will continue to do so from now on. A superscript $k$ denotes a function evaluation at $\vec{\boldsymbol \tau}^k$.

The system \eqref{eq-trs-newton} can be identified with the KKT conditions of a linear-quadratic problem. Indeed, one could equivalently try to obtain a step $\vec{\boldsymbol{\delta\tau}}$ by solving
\begin{subequations} \label{eq-trs-gen}
\begin{align}
\min \vec{\boldsymbol{\delta\tau}}^\top \nabla J^k + \frac{1}{2} \vec{\boldsymbol{\delta\tau}}^\top \nabla^2 J^k \vec{\boldsymbol{\delta\tau}}\\
\text{subject to} \quad  D \vec{\boldsymbol{\delta\tau}} + D \vec{\boldsymbol{\tau}}^k - f &= 0. \label{eq-trs-gen-com}
\end{align}
\end{subequations}

The linearity of the constraint \eqref{eq-tt-hb-fds-mp-com} allows us to implement some simplifications to the general numerical method. Given any $\vec{\boldsymbol \tau} \in \mathbb{R}^{2n_\mathrm{T}}$, we can find a projection $\vec{\boldsymbol \tau}^0$ onto the manifold $\lbrace \vec{\boldsymbol \tau} \in \mathbb{R}^{2n_\mathrm{T}}: D \vec{\boldsymbol \tau} = f \rbrace$ by setting
\begin{equation} \label{eq-trs-prj}
\vec{\boldsymbol \tau}^0 := \vec{\boldsymbol \tau} - A^{-1}D^\top(DA^{-1}D^\top)^{-1}(D\vec{\boldsymbol \tau} - f),
\end{equation}
where $A$ is a diagonal matrix that contains the triangle areas
\begin{equation*}
A := \mathrm{diag}(\vert T_1 \vert, \vert T_1 \vert, \dots, \vert T_{n_\mathrm{T}} \vert, \vert T_{n_\mathrm{T}} \vert).
\end{equation*}
The matrix $DA^{-1}D^\top$ is the usual discretisation of the Laplacian with piecewise linear Lagrange elements and $A^{-1}D^\top$ is a discrete analogue of the gradient.

With \eqref{eq-trs-prj}, the constraint \eqref{eq-trs-gen-com} simplifies to
\begin{equation}
D \vec{\boldsymbol{\delta\tau}} = 0 \label{eq-trs-smpl-com}
\end{equation}
and all iterates $\vec{\boldsymbol \tau}^k$, $k=0,1,2,...$ are automatically feasible with respect to \eqref{eq-tt-hb-fds-mp-com}.

Let $(\vec{\boldsymbol \tau}^*, \vec{y}^*) \in \mathbb{R}^{n_\mathrm{T}}\times \mathbb{R}^{n_\mathrm{I}}$ be an exact solution of the KKT conditions \eqref{eq-tt-kkt-fds}. Given an iterate $\vec{\boldsymbol \tau}^k \in \mathbb{R}^{n_\mathrm{T}} $, the overdetermined system \eqref{eq-tt-kkt1-fds} may not have a corresponding solution $\vec{y}^k \in \mathbb{R}^{n_\mathrm{I}}$. However, we can obtain a solution in the least-square sense by solving
\begin{equation}
\vec{y}^k := \mathrm{arg\,min}\; \vert \nabla J^k - D^\top \vec{y} \vert^2 = (DD^\top)^{-1}D\nabla J^k \label{eq-trs-yest}.
\end{equation}
Later on, we will use \eqref{eq-trs-yest} to get estimates of $\vec{y}^*$ from the iterates $\vec{\boldsymbol \tau}^k$.

\subsection{Trust region}

Trust-region methods impose an additional constraint on an optimisation problem by limiting the step size with a trust radius $\Delta_{\mathrm{TR}}^k > 0$. We apply the trust-region concept to the subproblem \eqref{eq-trs-gen} with \eqref{eq-trs-smpl-com}. Intuitively, one could say that we only trust the Taylor approximation of the objective in \eqref{eq-trs-gen} to be a sufficiently accurate representation of the exact objective within a ball of radius $\Delta_{\mathrm{TR}}^k$ around $\vec{\boldsymbol \tau}^k$.
\begin{subequations} \label{eq-trs-tan}
\begin{align}
\min \vec{\boldsymbol{\delta\tau}}^\top \nabla J^k + \frac{1}{2} \vec{\boldsymbol{\delta\tau}}^\top & \nabla^2 J^k \vec{\boldsymbol{\delta\tau}}\label{eq-trs-tan-obj}\\ 
\text{subject to} \quad  D \vec{\boldsymbol{\delta\tau}} &= 0 \label{eq-trs-tan-com}\\ 
\vert \vec{\boldsymbol{\delta\tau}} \vert &\leq \Delta_\mathrm{TR}^k. \label{eq-trs-tan-tr}
\end{align}
\end{subequations}
Through this amendment, we enforce a finite step length even if the Hessian is not positive definite. From \eqref{eq-tt-hess}, we observe that zero eigenvalues will occur whenever $\vert \vec{\boldsymbol \tau}_i \vert \leq \tau_0$ on a triangle $T_i$ ($i\in\lbrace 1,...,n_\mathrm{T} \rbrace$). Hence, we can only assume the Hessian $\nabla^2 J$ to be positive semidefinite, but not positive definite.

Generally, trust-region SQP methods require the solution of two subproblems, often referred to as the vertical (or normal) and the horizontal (or tangential) subproblems. With an affine constraint like \eqref{eq-trs-gen-com}, the set of admissible $\vec{\boldsymbol{\delta\tau}}$ that satisfy both this constraint and the trust-region constraint may be empty. Then, the purpose of the auxiliary vertical subproblem is to replace the equality constraint with a suitable relaxation. However, due to the initial projection that we suggest, \eqref{eq-trs-tan-com} and \eqref{eq-trs-tan-tr} are clearly not mutually exclusive ($\vec{\boldsymbol{\delta\tau}} = 0$, for instance, satisfies both constraints). As a result, for the problem under consideration, we may do without the vertical subproblem. Without any further modifications, the horizontal subproblem is given by \eqref{eq-trs-tan}.

\subsection{Solution of the horizontal subproblem}

We are looking for an efficient method that provides an approximate solution to the horizontal subproblem \eqref{eq-trs-tan}. While the classical Conjugate Gradient (CG) method assumes a positive definite Hessian, Steihaug \cite{Steihaug1983} and Toint \cite{Toint1981} developed an extension that remains valid for positive semidefinite or even indefinite matrices. Our implementation of the projected CG-Steihaug algorithm is an adaptation of the algorithms in \cite{Lalee1998,Nocedal2006}.

\begin{algorithm}[H]
\TitleOfAlgo{Projected CG-Steihaug (CGS)}
\KwData{$k\in\mathbb{N}_0$, $\mathtt{divtol} > 0, \mathtt{abstol} > 0, \mathtt{reltol} > 0,\Delta_\mathrm{TR}^k > 0$}
Set $\vec{\boldsymbol z}^0 = 0$, $\vec{\boldsymbol r}^0 = \nabla J^k$, $\vec{\boldsymbol g}^0 = \vec{\boldsymbol r}^0 - D^\top (DD^\top)^{-1}D\vec{\boldsymbol r}^0$, $\vec{\boldsymbol d}^0 = -\vec{\boldsymbol g}^0$\;
\lIf{$\sqrt{{\vec{\boldsymbol g}^{0\top}} \vec{\boldsymbol r}^0} < \mathtt{abstol}$}{\Return $\vec{\boldsymbol {\delta \tau}^k} = 0$}
\For{$j=0,1,2,...$}{
\If{$\vec{\boldsymbol d}^{j\top} \nabla^2 J^k \vec{\boldsymbol d}^{j} < \mathtt{divtol}$}{
\tcc{Small or zero curvature}
\tcc{Intersect step with trust region boundary}
Find $s>0$ such that $\vert \vec{\boldsymbol z}^{j} + s \vec{\boldsymbol d}^{j} \vert^2 = \Delta_\mathrm{TR}^2$\;
\While{Armijo condition \eqref{eq-trs-armijo} violated}{
Backtracking $s\to s/2$
}
\Return $\vec{\boldsymbol {\delta \tau}^+} = \vec{\boldsymbol z}^{j} + s \vec{\boldsymbol d}^{j}$\;
}
Compute $\vec{\boldsymbol z}^{+} = \vec{\boldsymbol z}^{j} + \left({\vec{\boldsymbol g}^{j\top}} \vec{\boldsymbol r}^j\right)/\left(\vec{\boldsymbol d}^{j\top} \nabla^2 J^k \vec{\boldsymbol d}^{j}\right) \vec{\boldsymbol d}^{j}$\;
\If{$\vert \vec{\boldsymbol z}^{+} \vert \geq \Delta_\mathrm{TR}^k$}{
\tcc{Step exceeds trust region}
\tcc{Intersect step with trust region boundary}
Find $s>0$ such that $\vert \vec{\boldsymbol z}^{j} + s \vec{\boldsymbol d}^{j} \vert^2 = \Delta_\mathrm{TR}^2$\;
\Return $\vec{\boldsymbol {\delta \tau}^+} = \vec{\boldsymbol z}^{j} + s \vec{\boldsymbol d}^{j}$\;
}
Set $\vec{\boldsymbol z}^{j+1} = \vec{\boldsymbol z}^{+}$, $\vec{\boldsymbol r}^{j+1} = \vec{\boldsymbol r}^{j} + \left({\vec{\boldsymbol g}^{j\top}} \vec{\boldsymbol r}^j\right)/\left(\vec{\boldsymbol d}^{j\top} \nabla^2 J^k \vec{\boldsymbol d}^{j}\right) \nabla^2 J^k \vec{\boldsymbol d}^{j}$, $\vec{\boldsymbol g}^{j+1} = \vec{\boldsymbol r}^{j+1} - D^\top (DD^\top)^{-1}D\vec{\boldsymbol r}^{j+1}$\;
\If{$\sqrt{{\vec{\boldsymbol g}^{j+1\top}} \vec{\boldsymbol r}^{j+1}} < \mathtt{reltol} \sqrt{{\vec{\boldsymbol g}^{0\top}} \vec{\boldsymbol r}^0}$}{
\tcc{Solution has converged}
\Return $\vec{\boldsymbol {\delta \tau}^+} = \vec{\boldsymbol z}^{j+1}$\;
}
Set $\beta^{j+1} = \left({\vec{\boldsymbol g}^{j+1\top}} \vec{\boldsymbol r}^{j+1}\right)/\left({\vec{\boldsymbol g}^{j\top}} \vec{\boldsymbol r}^{j}\right)$\;
Compute $\vec{\boldsymbol d}^{j+1} = -\vec{\boldsymbol g}^{j+1} + \beta^{j+1} \vec{\boldsymbol d}^{j}$\;
}
\end{algorithm}

The Armijo condition
\begin{equation}
Q(\vec{\boldsymbol z}^{j} + s \vec{\boldsymbol d}^{j}) \leq Q(\vec{\boldsymbol z}^{j}) + \gamma s \nabla Q(\vec{\boldsymbol z}^{j})^\top \boldsymbol d^j \label{eq-trs-armijo}
\end{equation}
with $\gamma = 1$e-2 guarantees sufficient decrease in the model objective $Q(\boldsymbol z) := \vec{\boldsymbol z}^\top \nabla J^k + \frac{1}{2} \vec{\boldsymbol z}^\top \nabla^2 J^k \vec{\boldsymbol z}$ even if very small curvature is encountered.

One can show that while the model objective \eqref{eq-trs-tan-obj} decreases monotonically with every iteration, the iterates strictly increase in their norms, $\vert \vec{\boldsymbol z}^{j+1} \vert > \vert \vec{\boldsymbol z}^{j} \vert$ \cite[p. 172]{Nocedal2006}. These properties justify the stopping criteria of the first two \textbf{if} statements within the \textbf{for} loop.

\subsection{Update of the trust radius}

To analyse the quality of a step $\vec{\boldsymbol {\delta \tau}}^+$ computed by Algorithm CGS, we compare the actual reduction in the value of the objective
\begin{align}
\mathtt{ared} &:= J^{k} - J^{k+1} \label{eq-trs-ared}
\intertext{with the reduction predicted by the quadratic model \eqref{eq-trs-tan-obj}}
\mathtt{pred} &:= -\vec{\boldsymbol{\delta\tau}}^{k\top} \nabla J^k - \frac{1}{2} \vec{\boldsymbol{\delta\tau}}^{k\top} \nabla^2 J^k \vec{\boldsymbol{\delta\tau}}^k. \label{eq-trs-pred}
\end{align}
If the actual reduction is sufficiently large, $\mathtt{ared} \geq \eta \; \mathtt{pred}$ with some parameter $\eta \in (0,1)$, then the step is accepted. A particularly good agreement of $\mathtt{ared}$ with $\mathtt{pred}$ indicates that we may even extend the trust region in the next step.

Otherwise, if $\mathtt{ared} < \eta \; \mathtt{pred}$, we infer that the discrepancy between the exact objective and its Taylor approximation is too large. Consequently, we reject the step computed in Algorithm CGS and shrink the trust region.

Following \cite{Lalee1998}, we update the trust radius as described below.

\begin{algorithm}[H]
\TitleOfAlgo{Update of the trust radius (UTR)}
\KwData{$k\in\mathbb{N}_0$, $\eta \in (0,1)$, $\Delta_\mathrm{TR}^k > 0$, $\Delta_\mathrm{TR}^\mathrm{max} > 0$}
Compute $\rho = \frac{\mathtt{ared}}{\mathtt{pred}}$\;
\uIf{$\rho\geq \eta$}{
\tcc{Step is accepted}
Set $\vec{\boldsymbol {\delta \tau}}^{k} = \vec{\boldsymbol {\delta \tau}}^+$\;
\tcc{Increase trust radius depending on model accuracy}
\uIf{$\rho\geq 0.9$}{
$\Delta_\mathrm{TR}^{k+1} =\min(\max(10\vert \vec{\boldsymbol {\delta \tau}}^{k} \vert ,\Delta_\mathrm{TR}^k),\Delta_\mathrm{TR}^\mathrm{max})$\;
}
\ElseIf{$\rho\geq 0.3$}{
$\Delta_\mathrm{TR}^{k+1} =\min(\max(2\vert \vec{\boldsymbol {\delta \tau}}^{k} \vert ,\Delta_\mathrm{TR}^k),\Delta_\mathrm{TR}^\mathrm{max})$\;
}
}
\Else{
\tcc{Step is rejected}
Set $\vec{\boldsymbol {\delta \tau}}^{k} = 0$\;
\tcc{Decrease trust radius}
$\Delta_\mathrm{TR}^{k+1} =\max(0.1,\min(0.5,\frac{1-\eta}{1-\rho}))\vert \vec{\boldsymbol {\delta \tau}}^{k} \vert$\;
}
\end{algorithm}

\subsection{Algorithm TRS}

All in all, the computation of a solution to \eqref{eq-tt-hb-fds-mp} follows the iterative scheme as outlined below.

\begin{algorithm}[H]
\TitleOfAlgo{Trust-region SQP (TRS)}
\KwData{$\mathtt{abstol}, \mathtt{reltol} > 0, \eta \in (0,1), \Delta_\mathrm{TR}^\mathrm{max} > 0, \Delta_\mathrm{TR}^0 \in (0,\Delta_\mathrm{TR}^\mathrm{max}] $}
Choose $\vec{\boldsymbol \tau}\in\mathbb{R}^{2n_T}$\;
Compute projection $\vec{\boldsymbol \tau}^0$ with \eqref{eq-trs-prj}\;
\For{$k = 0,1, 2, ...$}{
Evaluate $\nabla J^k$ and compute $\vec{y}^k$ with \eqref{eq-trs-yest}\;
\If{$\vert \nabla J^k - D^\top \vec{y}^k \vert_\infty \leq \mathtt{abstol}$ \KwSty{and} $\vert \vec{y}^k - \vec{y}^{k-1} \vert \leq \mathtt{reltol} \vert \vec{y}^k\vert$}{\Return $\vec{\boldsymbol \tau}^k, \vec{y}^k$}
Evaluate $\nabla^2 J^k$\;
Solve the tangential subproblem \eqref{eq-trs-tan} for $\vec{\boldsymbol{\delta\tau}}^+$ with Algorithm CGS\;
Compute $\vec{\boldsymbol\tau}^+ := \vec{\boldsymbol\tau}^{k}+\vec{\boldsymbol{\delta\tau}}^+$\;
Evaluate $J^k$ and $J^{+}$\;
Compute $\mathtt{ared}$ and $\mathtt{pred}$ from \eqref{eq-trs-ared} and \eqref{eq-trs-pred}\;
Use Algorithm UTR to accept or reject $\vec{\boldsymbol{\delta\tau}}^+$\;
}
\end{algorithm}

We note that a stopping criterion which measures the relative difference between subsequent iterates of the stress is not appropriate, since in general the stress is not uniquely determined.

\section{Numerical experiments}
\label{sec-experiments}

We will now present numerical results to demonstrate the performance of a \textsf{MATLAB} implementation of Algorithm TRS.

In all examples, we initialise the algorithm with $\vec{\boldsymbol \tau} = 0$, i.e. the initial iterate is given by $\vec{\boldsymbol \tau}^0 = -A^{-1}D^\top (DA^{-1}D^\top)^{-1}f$ according to \eqref{eq-trs-prj}. Furthermore, we set \texttt{abstol} = \texttt{reltol} = 1e-4, \texttt{divtol} = 1e-10, $\Delta_\mathrm{TR}^0$ = 10, $\Delta_\mathrm{TR}^\mathrm{max}$ = 1e5 and $\eta$ = 1e-1. For the problem parameters we assume $\mu$ = $\kappa$ = 1 and $f = 1$.

Numerical solutions for Bingham and Herschel-Bulkley flows through various duct geometries are available in the literature \cite{Moyers2004,Huilgol2005}. These were computed with Algorithm ALG2, an Uzawa type method based on an augmented Lagrangian formulation of the nonsmooth minimisation problem \eqref{eq-dl-mp}. As \cite{Huilgol2005} already contains a discussion of the results in \cite{Moyers2004}, we will use the data provided in \cite{Huilgol2005} as reference solutions.

In order to compare the runtimes of Algorithms TRS and ALG2, we also provide results which we obtained with our own \textsf{MATLAB} implementation of Algorithm ALG2. Here we use the discrete analogue of the augmented Lagrangian
\begin{equation} \label{eq-alg}
L_r(y,\boldsymbol q,\boldsymbol \tau) = c(y) + j(y) - (f,y)_\alpha - (\boldsymbol \tau, \boldsymbol q - \nabla y)_\alpha + \frac{r}{2} \int\limits_\Omega \vert \nabla y - \boldsymbol q \vert^2 \; \mathrm{d}x,
\end{equation}
discretised with the same finite elements as for Algorithm TRS. We remark that for $\alpha < 2$ it is not obvious whether the penalty term in \eqref{eq-alg} is well-defined since $\nabla y$ and $\boldsymbol q$ do not necessarily belong to $\boldsymbol L^2 (\Omega)$. This problem dissolves in the discrete setting. Alternatively one might consider a non-quadratic penalty term, e.g. by replacing the exponent $2$ with $\alpha$. In that case, however, the penalty term will not give rise to a Laplacian in the optimality conditions, but rather a quasilinear elliptic operator similar to an $\alpha$-Laplacian. Consequently, ALG2 would require the solution of an additional nonlinear equation in every single iteration.

We set $r$ = 10, which by trial and error appears to be an optimal choice to minimise the number of iterations. To ensure consistency with the set up of Algorithm TRS, we terminate Algorithm ALG2 as soon as the residual of the first-order optimality conditions measured in the $\infty$-norm no longer exceeds \texttt{abstol}. We also require subsequent iterates of the discretised velocity $\vec{y}$ and the discretised rate of strain $\vec{\boldsymbol q}$ to have a relative difference of at most \texttt{reltol}. For Newton's method, which is required to find a new iterate for $\vec{\boldsymbol q}$ if $\alpha \neq 2$ (cf \cite{Huilgol2005}), we also use \texttt{abstol} and \texttt{reltol} as absolute and relative tolerances.

Our meshes are generated with the \textsf{MATLAB} built-in Delaunay triangulation. The programs are executed with \textsf{MATLAB} R2013a 64-bit on a laptop with an \textsf{Intel\textsuperscript\textregistered Core\texttrademark i5} CPU 4x2.53 GHz.

\subsection{Flow through a cylindrical pipe}

\begin{figure}[p]
\centering
\includegraphics[width=0.75\textwidth]{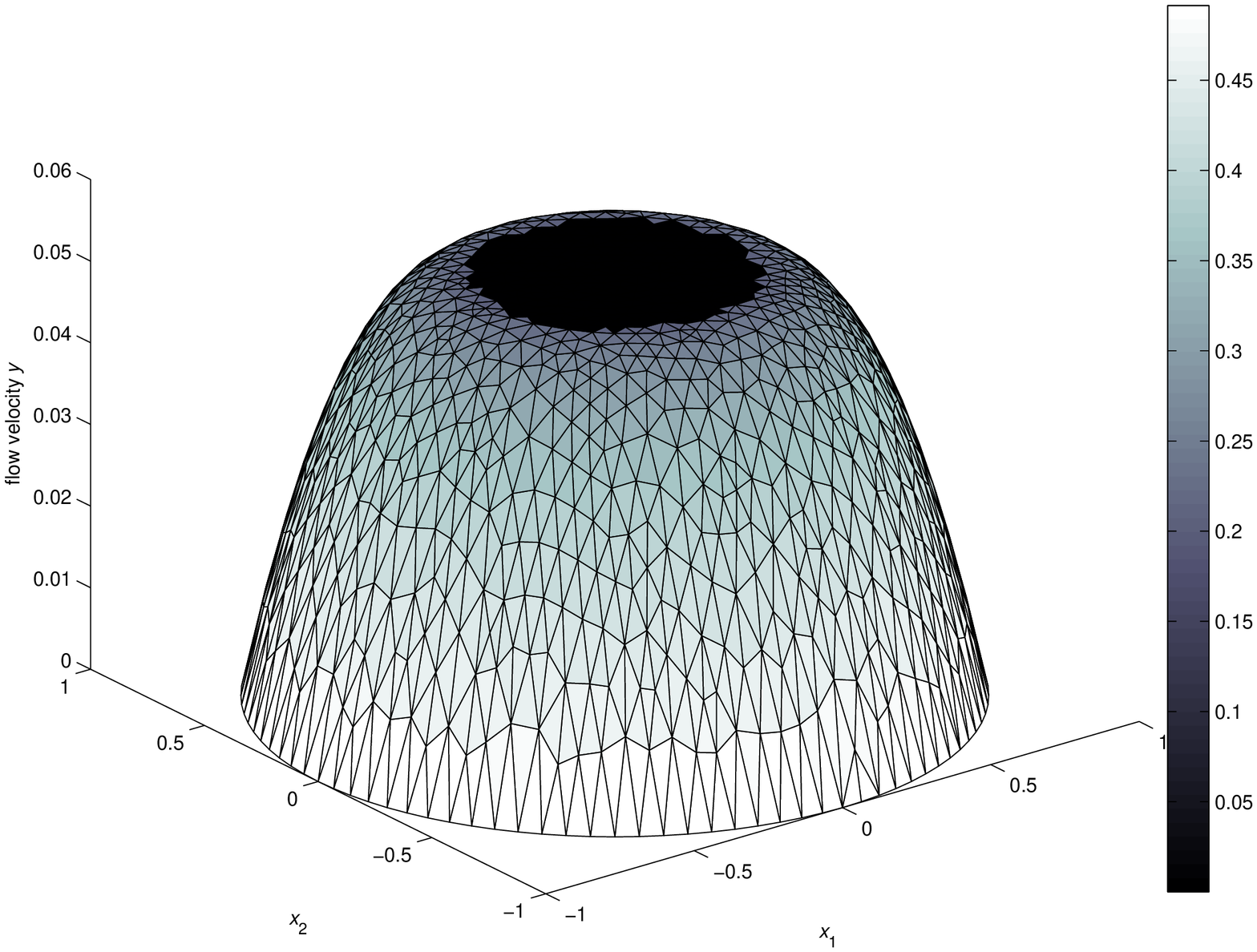}
\caption{Typical velocity profile for viscoplastic duct flow, computed with Algorithm TRS. Here, $\alpha = 1.75$ and $\tau_0 = 0.2$. On triangles $T_i$ that are coloured black, $\vert \vec{\boldsymbol \tau}_i \vert \leq \tau_0$. Otherwise, the colours represent $\vert \vec{\boldsymbol \tau}_i \vert$.}
\label{fig-cyl-surface}
\end{figure}

We consider the classical test problem of flow through a cylindrical pipe with radius 1 and homogeneous Dirichlet boundary conditions. A typical solution is depicted in Figure \ref{fig-cyl-surface}. The graph illustrates very well the existence of a plastic region with stagnant flow which is separated from a region with regular viscous flow.

For this quasi one-dimensional problem, an analytical solution is known. With $R := \vert x \vert$, $R_0 := 2\tau_0$ and $\beta := 1/(\alpha - 1)$ it reads \cite{Froishteter1980}
\begin{equation}
y(R) =
\begin{cases}
\frac{1}{2^\beta(1+\beta)}(1-R_0)^{1+\beta} & 0 \leq R \leq R_0\\
\frac{1}{2^\beta(1+\beta)}\left((1-R_0)^{1+\beta}-(R-R_0)^{1+\beta}\right) & R_0 < R \leq 1
\end{cases}.
\end{equation}

Huilgol and You \cite{Huilgol2005} provide results which they computed with their implementation of Algorithm ALG2 for this problem. They consider the cases of $\alpha = 2$, $\alpha = 1.75$ and $\alpha = 1.5$ on three different grids and with $\tau_0 = 0.1$ or $\tau_0 = 0.2$. The authors write $n$ for $\alpha - 1$, $Od$ for $2\tau_0$ and $u/U$ for $\vec{y}$.

While Huilgol and You report that a ``convergence tolerance $\epsilon = 10^{-4}$ is kept fixed'' \cite[p 130]{Huilgol2005}, it unfortunately remains unclear what this tolerance refers to. Their numerical results suggest that their stopping criterion only measures a relative difference between subsequent iterates of the flow velocity. In this case, their implementation of ALG2 would neglect to check how well the first order necessary optimality conditions of the minimisation problem are satisfied.

In Tables \ref{tb-alg-cyl-bi} to \ref{tb-alg-cyl-hb2} we present some typical results of Algorithm TRS and compare them with those of our implementation of Algorithm ALG2 as well as the analytical solution. To verify that Algorithms TRS and ALG2 compute accurate solutions, we compute a relative error $\vert \vec{y}- \vec{y}_a \vert/\vert \vec{y}_a \vert$ between the converged solution $\vec{y}$ and the analytical solution $\vec{y}_a$ evaluated on the vertices of the mesh. Furthermore, we provide the number of iterations required to satisfy each stopping criterion, the CPU times needed for each algorithm to terminate and the speedup as the ratio of the times for ALG2 and TRS.

\begin{table}
\begin{tabular}{c|c|cc|cc|cc|c}
\hline
& & \multicolumn{2}{c|}{Relative error} & \multicolumn{2}{c|}{Iterations} & \multicolumn{2}{c|}{CPU time (s)} & \\
$\tau_0$ & $n_\mathrm{N}$ & TRS & ALG2 & TRS & ALG2 & TRS & ALG2 & Speedup\\
\hline
0.1 & 559 & 1.48e-3 & 2.27e-3 & 3 & 73 & 0.06 & 0.36 & 6.0\\
0.1 & 1129 & 6.35e-4 & 1.42e-3 & 2 & 73 & 0.06 & 0.66 & 11\\
0.1 & 2169 & 3.40e-4 & 1.24e-3 & 1 & 73 & 0.07 & 1.19 & 17\\
0.2 & 559 & 2.19e-3 & 3.01e-3 & 58 & 73 & 0.49 & 0.40 & 0.82\\
0.2 & 1129 & 8.97e-4 & 1.75e-3 & 10 & 73 & 0.16 & 0.66 & 4.1\\
0.2 & 2169 & 5.30e-4 & 1.41e-3 & 2 & 73 & 0.09 & 1.18 & 13\\
\hline
\end{tabular}
\caption{Bingham flow through the cylindrical pipe, $\alpha = 2$, using $n_\mathrm{N}$ triangle nodes.}
\label{tb-alg-cyl-bi}
\end{table}

\begin{table}
\begin{tabular}{c|c|cc|cc|cc|c}
\hline
& & \multicolumn{2}{c|}{Relative error} & \multicolumn{2}{c|}{Iterations} & \multicolumn{2}{c|}{CPU time (s)} & \\
$\tau_0$ & $n_\mathrm{N}$ & TRS & ALG2 & TRS & ALG2 & TRS & ALG2 & Speedup\\
\hline
0.1 & 559 & 1.97e-3 & 2.67e-3 & 11 & 67 & 0.14 & 0.76 & 5.4\\
0.1 & 1129 & 8.79e-4 & 1.55e-3 & 2 & 67 & 0.07 & 1.36 & 19\\
0.1 & 2169 & 4.56e-4 & 1.23e-3 & 1 & 67 & 0.08 & 2.55 & 32\\
0.2 & 559 & 3.03e-3 & 3.68e-3 & 51 & 62 & 0.46 & 0.70 & 1.5\\
0.2 & 1129 & 1.33e-4 & 1.98e-3 & 15 & 62 & 0.25 & 1.24 & 5.0\\
0.2 & 2169 & 7.28e-4 & 1.41e-3 & 20 & 62 & 0.61 & 2.54 & 4.2\\
\hline
\end{tabular}
\caption{Herschel-Bulkley flow through the cylindrical pipe, $\alpha = 1.75$.}
\label{tb-alg-cyl-hb1}
\end{table}

\begin{table}
\begin{tabular}{c|c|cc|cc|cc|c}
\hline
& & \multicolumn{2}{c|}{Relative error} & \multicolumn{2}{c|}{Iterations} & \multicolumn{2}{c|}{CPU time (s)} & \\
$\tau_0$ & $n_\mathrm{N}$ & TRS & ALG2 & TRS & ALG2 & TRS & ALG2 & Speedup\\
\hline
0.1 & 559 & 3.38e-3 & 3.84e-3 & 20 & 54 & 0.20 & 1.82 & 9.1\\
0.1 & 1129 & 1.50e-3 & 2.01e-3 & 3 & 54 & 0.08 & 9.87 & 120\\
0.1 & 2169 & 7.87e-4 & 1.32e-3 & 3 & 54 & 0.12 & 23.02 & 190\\
0.2 & 559 & 5.68e-3 & 6.05e-3 & 37 & 43 & 0.34 & 4.88 & 14\\
0.2 & 1129 & 2.69e-3 & 3.06e-3 & 11 & 43 & 0.17 & 13.41 & 79\\
0.2 & 2169 & 1.46e-3 & 1.81e-3 & 12 & 43 & 0.31 & 27.89 & 90\\
\hline
\end{tabular}
\caption{Herschel-Bulkley flow through the cylindrical pipe, $\alpha = 1.5$.}
\label{tb-alg-cyl-hb2}
\end{table}

Looking at the relative errors in Tables \ref{tb-alg-cyl-bi} to \ref{tb-alg-cyl-hb2}, we conclude that both algorithms compute accurate approximations to the analytical solution. Errors tend to decrease as the mesh is refined. This behaviour indicates convergence of the numerical solutions to the exact solution as $h \to 0$. Our results are also consistent with the data in \cite{Huilgol2005}.

Moreover, it takes consistently fewer iterations for Algorithm TRS to converge compared to Algorithm ALG2. With only one single exception for Bingham flow, the runtime for Algorithm TRS amounts to just a small fraction of the time Algorithm ALG2 takes to converge. This enormous difference in performance of the two methods increases even further in favour of the trust-region SQP algorithm when $\alpha$ is decreased. Overall, Algorithm TRS is hardly affected by changes in the power-law exponent $\alpha$, neither in terms of iterations, nor in terms of CPU times. In sharp contrast, the number of iterations decreases slightly for Algorithm ALG2 as $\alpha$ decreases, while the computational expenditure per iteration increases significantly. This is due to the fact that in the Bingham case, the auxiliary variable $\boldsymbol q$ is obtained through a simple function evaluation. With general Herschel-Bulkley fluids, however, Newton's method is required to solve a nonlinear equation for $\boldsymbol q$ in the region of the domain where the fluid has yielded. The more $\alpha$ deviates from the value 2, the more iterations of Newton's method are needed in every single step of Algorithm ALG2 to find the next iterate of $\boldsymbol q$.

An interesting feature of the augmented Lagrangian method is its remarkable mesh independence. For all combinations of the yield stress $\tau_0$ and the exponent $\alpha$ chosen here, the number of iterations is exactly the same for all three grids. The trust-region SQP method exhibits a trend of a decreasing number of iterations as the mesh is refined, and in some cases even terminates after the first iteration. With a simple calculation, one verifies for the cylindrical duct that one solution of the continuous problem \eqref{eq-tt-mp} is given by $\boldsymbol \tau := -\nabla(-\Delta)^{-1}f$. Algorithm TRS is initialised with $\vec{\boldsymbol \tau} = A^{-1}D^\top (DA^{-1}D^\top)^{-1}f$, which corresponds to the finite element discretisation of this solution. Hence, with $h$ sufficiently small, this initial projection already solves the discrete optimality conditions \eqref{eq-tt-kkt-fds} and Algorithm TRS terminates.

Overall we conclude that both methods are suited to compute accurate approximations to the cylindrical duct flow. However, Algorithm TRS generally converges several times faster than Algorithm ALG2.

\subsection{Flow through a kiwi shaped duct}

One may argue that the advantageous performance of Algorithm TRS could, at least partially, be an artefact of the symmetry of the cylindrical duct. It is this symmetry that allows the problem to be reduced to one dimension and to solve it analytically. In order to eliminate such effects, we now consider a more complicated geometry with a non-symmetric domain.

\begin{figure}[p]
\centering
\includegraphics[width=0.5\textwidth]{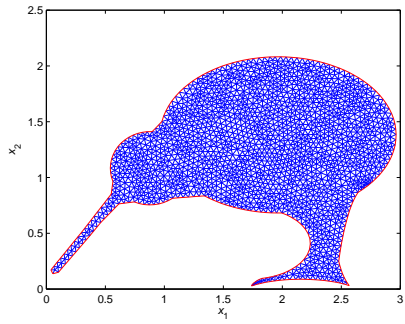}
\caption{Geometry of the kiwi shaped pipe with a sample mesh ($n_\mathrm{N} = 2101$).}
\label{fig-kiwimesh}
\end{figure}

\begin{figure}[p]
\centering
\includegraphics[width=0.75\textwidth]{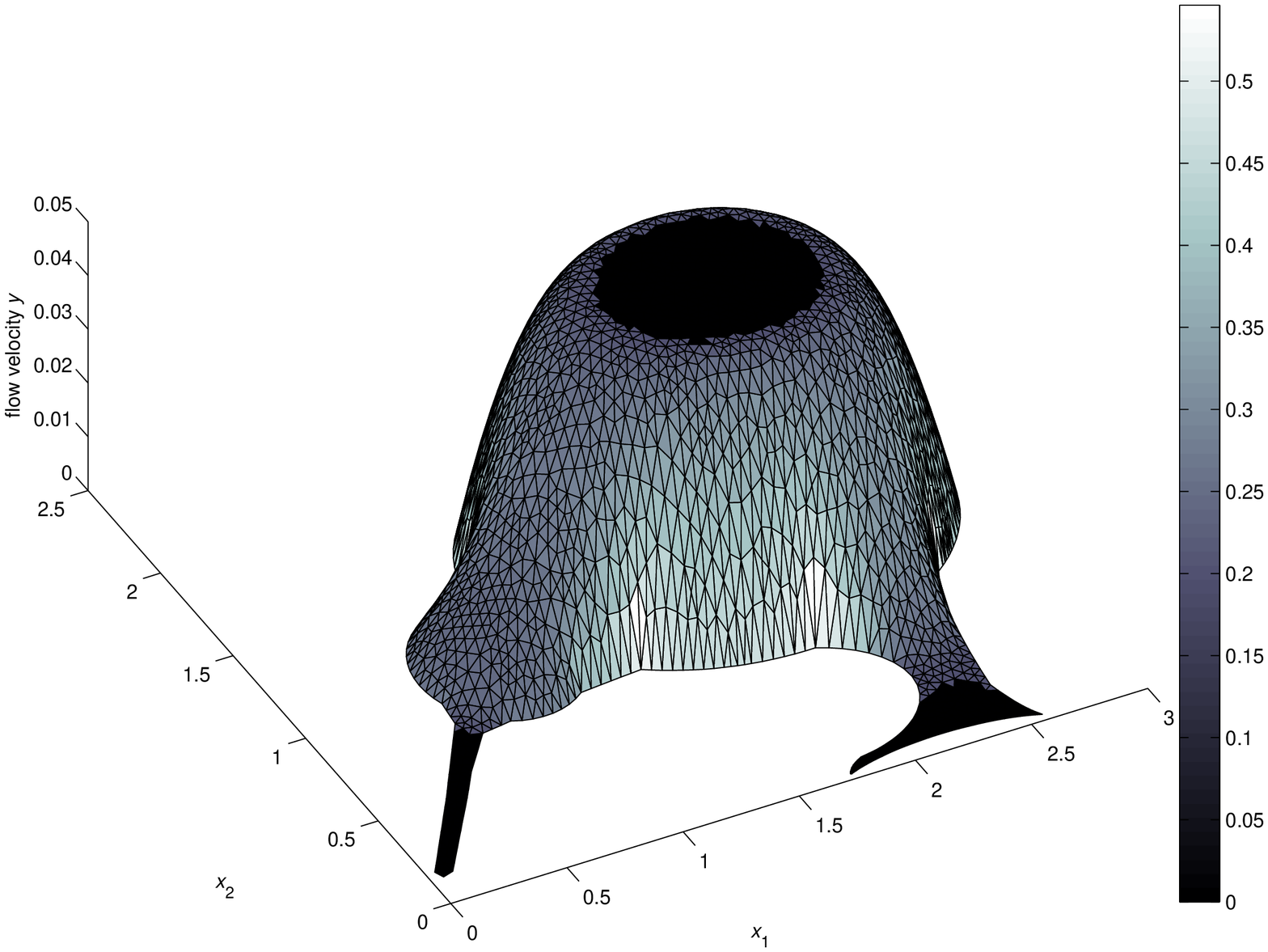}
\caption{TRS approximation of Bingham flow through the kiwi shaped pipe, $\alpha = 2$, $\tau_0 = 0.2$, $n_\mathrm{N} = 2101$.  As in Figure \protect\ref{fig-cyl-surface}, the colours indicate locations of stagnant flow and represent the stress magnitude in yielded regions.}
\label{fig-kiwi-trs}
\end{figure}

As shown in Figure \ref{fig-kiwimesh}, we now investigate Herschel-Bulkley flow through a pipe with the shape of a kiwi. Figure \ref{fig-kiwi-trs} depicts a typical solution of viscoplastic flow through a pipe with this cross-section and no-slip boundary conditions. In the thin regions of the domain around the beak and the foot, the fluid remains stuck and does not flow at all. Similarly, another black area in the centre of the body indicates how a column of fluid moves rigidly through the duct.

\begin{figure}[p]
\centering
\includegraphics[width=0.75\textwidth]{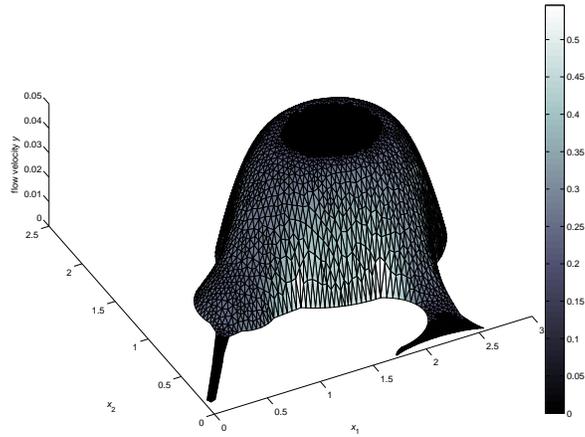}
\caption{ALG2 approximation of Bingham flow through the kiwi shaped pipe, $\alpha = 2$, $\tau_0 = 0.2$, $n_\mathrm{N} = 2101$.}
\label{fig-kiwi-alg}
\end{figure}

\begin{figure}[p]
\centering
\includegraphics[width=0.75\textwidth]{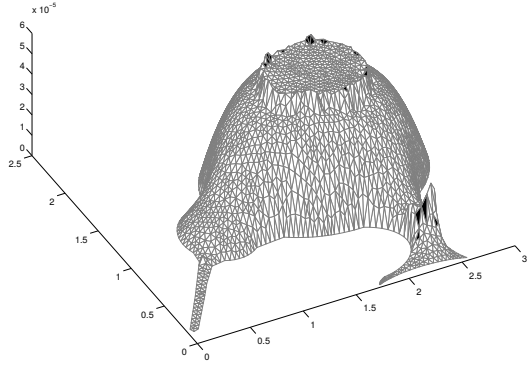}
\caption{Absolute error between the approximations of Figures \protect\ref{fig-kiwi-trs} and \protect\ref{fig-kiwi-alg}. Black triangles are identified as unyielded by Algorithm TRS, but not by Algorithm ALG2. White triangles are classified equally.}
\label{fig-kiwi-err}
\end{figure}

For the particular choice of parameters in Figures \ref{fig-kiwi-trs} and \ref{fig-kiwi-alg}, both velocity profiles appear to be identical by visual inspection. There is a slight discrepancy in the approximation of the solid regions: in the graph for Algorithm TRS, single triangles are coloured black unlike in the approximate solution generated by Algorithm ALG2. This insignificant difference is also illustrated in Figure \ref{fig-kiwi-err}. The surface in the latter graph represents the absolute error $\vert \vec{y}_\mathrm{TRS} - \vec{y}_\mathrm{ALG2} \vert$, the difference in the solutions obtained by the two methods. As expected, the difference reaches local highs near interfaces between viscous and plastic regions. In those plastic regions, the functional on which the augmented Lagrangian method is based is nonsmotth. Near the same interfaces, the gradients and Hessians which occur in the trust-region SQP algorithm become singular. As a result, the approximate solutions differ most in their prediction of the plug flow velocity.

In Tables \ref{tb-alg-kwi-bi} to \ref{tb-alg-kwi-hb2}, we present similar data for the simulations like before in Tables \ref{tb-alg-cyl-bi} to \ref{tb-alg-cyl-hb2}. Since an analytical solution is unavailable, we compute a relative difference $\vert \vec{y}_\mathrm{TRS} - \vec{y}_\mathrm{ALG2}\vert/\vert \vec{y}_\mathrm{ALG2}\vert$ instead of the relative error.

Even though in none of the examples Algorithm TRS terminates after the first iteration, its clearly superior performance persists even in this more complex geometry. The augmented Lagrangian method still exhibits mesh independence to a great extent. The number of iterations for the trust-region SQP method tends to decrease when the mesh is refined or when $\tau_0$ is decreased. No obvious correlation can be found between iterations of Algorithm TRS and the exponent $\alpha$. 

Again, there is only one exception to the overall trend that the computation with Algorithm TRS terminates significantly faster than Algorithm ALG2. With $\alpha = 1.75$ and even more with $\alpha = 1.5$, the computational cost of the augmented Lagrangian method becomes excessively high. Comparing the speedup factors for the two geometries, it appears that the symmetry of the problem is exploited more efficiently in the trust-region SQP method. However, even without any symmetry in the domain, we still achieve major gains in the computational performance by using our new method TRS.

\begin{table}
\begin{tabular}{c|c|c|cc|cc|c}
\hline
& & & \multicolumn{2}{c|}{Iterations} & \multicolumn{2}{c|}{CPU time (s)} & \\
$\tau_0$ & $n_\mathrm{N}$ & Relative difference & TRS & ALG2 & TRS & ALG2 & Speedup\\
\hline
0.1 & 1125 & 9.67e-4 & 14 & 73 & 0.20 & 0.53 & 2.7\\
0.1 & 2101 & 9.67e-4 & 13 & 73 & 0.29 & 0.96 & 3.3\\
0.1 & 4308 & 9.69e-4 & 8 & 73 & 0.34 & 2.16 & 6.4\\
0.2 & 1125 & 1.08e-3 & 80 & 73 & 1.00 & 0.55 & 0.55\\
0.2 & 2101 & 1.33e-3 & 27 & 73 & 0.50 & 0.94 & 1.9\\
0.2 & 4308 & 1.19e-3 & 37 & 73 & 1.29 & 2.15 & 1.7\\
\hline
\end{tabular}
\caption{Bingham flow through the kiwi pipe, $\alpha = 2$.}
\label{tb-alg-kwi-bi}
\end{table}

\begin{table}
\begin{tabular}{c|c|c|cc|cc|c}
\hline
& & & \multicolumn{2}{c|}{Iterations} & \multicolumn{2}{c|}{CPU time (s)} & \\
$\tau_0$ & $n_\mathrm{N}$ & Relative difference & TRS & ALG2 & TRS & ALG2 & Speedup\\
\hline
0.1 & 1125 & 7.28e-4 & 9 & 64 & 0.16 & 1.10 & 6.9\\
0.1 & 2101 & 7.38e-4 & 20 & 64 & 0.45 & 1.81 & 4.0\\
0.1 & 4308 & 7.42e-4 & 14 & 64 & 0.60 & 4.01 & 6.7\\
0.2 & 1125 & 9.09e-4 & 33 & 58 & 0.49 & 0.92 & 1.9\\
0.2 & 2101 & 7.92e-4 & 42 & 58 & 0.75 & 0.15 & 2.0\\
0.2 & 4308 & 1.46e-3 & 19 & 58 & 0.68 & 3.34 & 4.9\\
\hline
\end{tabular}
\caption{Herschel-Bulkley flow through the kiwi pipe, $\alpha = 1.75$.}
\label{tb-alg-kwi-hb1}
\end{table}

\begin{table}
\begin{tabular}{c|c|c|cc|cc|c}
\hline
& & & \multicolumn{2}{c|}{Iterations} & \multicolumn{2}{c|}{CPU time (s)} & \\
$\tau_0$ & $n_\mathrm{N}$ & Relative difference & TRS & ALG2 & TRS & ALG2 & Speedup\\
\hline
0.1 & 1125 & 3.63e-4 & 31 & 47 & 0.46 & 11.95 & 26\\
0.1 & 2101 & 4.13e-4 & 35 & 47 & 0.67 & 27.47 & 41\\
0.1 & 4308 & 6.47e-4 & 17 & 54 & 0.66 & 92.89 & 140\\
0.2 & 1125 & 3.21e-4 & 62 & 72 & 0.79 & 13.55 & 17\\
0.2 & 2101 & 5.46e-4 & 48 & 56 & 0.84 & 22.95 & 27\\
0.2 & 4308 & 1.10e-3 & 27 & 64 & 1.02 & 116.62 & 114\\
\hline
\end{tabular}
\caption{Herschel-Bulkley flow through the kiwi pipe, $\alpha = 1.5$.}
\label{tb-alg-kwi-hb2}
\end{table}

\section{Outlook and Conclusions}
\label{sec-conclusions}

So far, we have successfully developed a new approach to model and simulate viscoplastic flow. Analytically, the proposed dual formulation is fully equivalent to the conventional definition in terms of a variational inequality, or nonsmooth minimisation problem.

From a numerical perspective, however, the feature of higher regularity of the problem is very desirable. With Algorithm TRS, we have designed a powerful alternative to the well-known Algorithm ALG2. The latter method suffers from impaired performance due to the lack of continuous derivatives, the inability to exploit extra regularity in the case of shear-thinning fluids and the need to solve potentially large nonlinear systems repeatedly. Algorithm TRS, however, displays superior performance in particular for shear-thinning Herschel-Bulkley fluids. For our numerical examples, we generally observe large speedup factors, depending on the flow domain.

The extension of Algorithm TRS to time-dependent or fully two-dimensional flow is straightforward. While this will affect the explicit operators and function spaces involved, the overall structure of the problem remains preserved.

As with the augmented Lagrangian approach, the governing equations of viscoplastic flow no longer correspond to the optimality conditions of a minimisation problem when certain additional terms are considered in the model. This includes for instance finite slip (Robin) boundary conditions, stick-slip boundary conditions \cite{Roquet2008} or convection. The key question for a generalisation of Algorithm TRS to include those features will be, how the trust radius can be updated when no objective is available. In analogy to our successful application of Algorithm TRS to flows of nonlinear shear-thinning fluids, the trust-region SQP method requires no costly inner loop to handle further nonlinearites. Hence, also for more complex flow problems, we expect significantly increased performance compared to the augmented Lagrangian method.

\bibliographystyle{elsarticle-num} 
\bibliography{bibliography}

\end{document}